\documentclass[12pt]{amsart}
\usepackage{wrapfig,epsfig,epsf,graphicx,amsmath,amssymb}

\newcommand{\Ass}{A_{\rm ss}}
\newcommand{\Ahatss}{\hat{A}_{\rm ss}}
\newcommand{\B}{\mathcal B}
\newcommand{\e}{\varepsilon}
\newcommand{\h}{h_{\rm ss}}

\newcommand{\hhat}{\hat{h}_{\rm ss}}

\newcommand{\E}{\mathcal E}
\newcommand{\f}{\mathcal F}
\newcommand{\g}{\mathcal G}

\newcommand{\LL}{{\mathcal L}}

\newcommand{\Pss}{P_{\! \rm ss}}

\newcommand{\R}{\mathbb R}
\newcommand{\T}{\mathbb T}
\renewcommand{\d}{\prime} 
\newcommand{\dd}{{\prime \prime}} 
\newcommand{\ddd}{{\prime \prime \prime}}

\renewcommand{\hbar}{\overline{h}}
\newcommand{\full}[2]{\frac{d #1}{d #2}}

\newtheorem{theorem}{Theorem}

\newtheorem{lemma}[theorem]{Lemma}
\newtheorem{proposition}[theorem]{Proposition}

\newtheorem{definition}{Definition}

\renewcommand{\part}[2]{\frac{\partial #1}{\partial #2}}

\begin{document}
\title[]
{Energy levels of steady states for thin film type equations}
\author[]
{R. S. Laugesen\thanks{laugesen@math.uiuc.edu} and
M. C. Pugh\thanks{mpugh@math.upenn.edu}}
\address{Department of Mathematics, University of Illinois, Urbana, IL
61801}
\address{Department of Mathematics, University of Pennsylvania,
 Philadelphia, PA 19104}
\date \today

\begin{abstract}
We study the phase space of the evolution equation
\[
h_t = -(f(h) h_{xxx})_x - (g(h) h_x)_x 
\]
by means of a dissipated energy (a Liapunov function).  Here $h(x,t)
\geq 0$, and at $h=0$ the coefficient functions $f>0$ and $g$ can
either degenerate to $0$, or blow up to $\infty$, or tend to a nonzero
constant.

We first show all positive periodic steady states are `energy
unstable' fixed points for the evolution (meaning the energy decreases
under some zero--mean perturbation) if $(g/f)^\dd \geq 0$ or if the
perturbations are allowed to have period longer than that of the
steady state.

For power law coefficients ($f(y) = y^n$ and $g(y) = \B y^m$ for some
$\B > 0$) we analytically determine the relative energy levels of
distinct steady states.  For example, with $m-n \in [1,2)$ and for
suitable choices of the period and mean value, we find three
fundamentally different steady states. The first is a constant steady
state that is nonlinearly stable and is a local minimum of the energy.
The second is a positive periodic steady state that is linearly
unstable and has higher energy than the constant steady state; it is a
saddle point. The third is a periodic collection of `droplet'
(compactly supported) steady states having lower energy than either
the positive steady state or the constant one.  Since the energy must
decrease along every orbit, these results significantly constrain the
dynamics of the evolution equation.

Our results suggest that heteroclinic connections could exist between
certain of the steady states, for example from the periodic steady
state to the droplet one. In a companion article we perform numerical
simulations to confirm their existence.
\end{abstract}

\maketitle
\baselineskip = 18pt
\section{Introduction}
\label{introduction}

We study the evolution equation
\begin{equation}  \label{evolve}
        h_t = - (f(h) h_{xxx})_x - (g(h)h_x)_x .
\end{equation}
This is the one dimensional version of $h_t = -
\nabla \cdot( f(h) \nabla \Delta h) - \nabla \cdot (g(h) \nabla h)$,
which has been used to model the dynamics of a thin film of viscous
liquid.  The air/liquid interface is at height $z = h(x,y,t) \geq 0$ and
the
liquid/solid interface is at $z=0$. The one dimensional equation (\ref{evolve}) applies if the liquid film is uniform in the $y$ direction. 

The fourth order term in the equation reflects surface tension effects
and the second order term can reflect gravity, van der Waals
interactions, thermocapillary effects or the geometry of the solid
substrate, for example.  Typically $f(h)=h^3 + \beta h^p$ where $0 < p
< 3, \beta \geq 0$, and $g(h) \sim \pm h^m$ as $h \rightarrow 0$,
where $m \in \R$. In certain applications $g(h)$ changes sign at some
positive $h$. We refer to \cite{myers,oron} for reviews of the
physical and modeling literature.

The extensively studied Cahn--Hilliard equation \cite{CH,vdW} also has the
form (\ref{evolve}), with $f \equiv 1$ and $g(h)=1-3h^2$.
See \cite{ABF,BF90,BF00,GNC99} for further references on the
Cahn--Hilliard equation.

Equations like (\ref{evolve}) are of mathematical as well as physical
interest: for example, Bertozzi and Pugh \cite{BPLW} conjectured that
blow-up ($||h(\cdot,t)||_\infty \to \infty$) is possible in some cases
({\it e.g.} if $f(h)=h^n, g(h)=h^m$ with $m > n+2$), and they have
proved
\cite{BPFS} {\it finite time} blow-up for $f(h) = h$ and $g(h) = h^m$
when $m \geq 3$. In \cite[\S8]{LP1} we related the steady states and
some of their properties to this blow-up conjecture.

\subsection*{Background and goals.} 
In \cite{LP2} we proved linear stability and instability results for
the positive {\em periodic} steady states of (\ref{evolve}).
Periodicity should not be regarded as a constraint, since if $f,g>0$
then positive bounded steady states must be periodic or constant, by
\cite[Theorem~B.1]{LP1}. And periodic steady states do exist for many
equations of type (\ref{evolve}), by the methods of \cite[\S2.2]{LP1}
or \cite{GNC95,GNC95a}, for example.  Compactly supported `droplet'
steady states only exist, though, if $g/f$ satisfies additional
constraints \cite[\S2.2]{LP1}, and can have relatively low regularity
at the contact line.

In this paper we concentrate mostly on positive periodic steady
states, and on droplet steady states with zero contact angles. Our
main investigative tool is the energy
\[
\E(h(\cdot,t)) = \int_0^X \left[ \frac{1}{2} h_x(x,t)^2 - H(h(x,t))
\right] dx ;
\]
here $h(x,t)$ is a smooth solution of (\ref{evolve}) that is
$X$-periodic in $x$, and $H(y)$ satisfies $H^\dd = g/f$. This energy
is strictly dissipated: $(d/dt) \E(h(\cdot,t)) \leq 0$ with equality
if and only if $h$ is a steady state (see
\S\ref{instability_defn}). Thus the energy is a Liapunov function for
the evolution.

We address two questions about the energy landscape of the evolution (\ref{evolve}).
\begin{enumerate}
\item Which steady states are local minima of the energy? Which are
saddle points?
\item 
Among steady states having the same period and the same area ({\it
i.e.} fluid volume), which 
has the lowest energy?
\end{enumerate}
Answering these questions
will help clarify the phase portrait of the evolution. For example, if
all small zero--mean perturbations of a steady state can be shown to
raise the energy, then the steady state might be asymptotically
stable: it might be that every smooth solution starting from near the
steady state relaxes back to the steady state as $t \to \infty$.  But
if some zero--mean perturbation {\em decreases} the energy, then
asymptotic stability definitely fails.

Our requirement that the perturbations have zero mean
seems reasonable from a physical standpoint, because it
corresponds to a disturbance of the fluid that alters the profile
without adding additional fluid. Mathematically it is reasonable
because the evolution equation (\ref{evolve}) preserves volume for
spatially periodic solutions: $\int h(x,t) \, dx = \int h(x,0) \, dx$
for all time $t$.  Thus zero--mean perturbations allow the possibility
of relaxation back to the original steady state, while nonzero--mean
perturbations do not.
 
\subsection*{A sketch of definitions and results.} 
Take $X > 0$.  If there is an $X$-periodic zero--mean perturbation $v$
such that $\E(\h+ \varepsilon v) < \E(\h)$ for all small $\epsilon>0$,
then we call the steady state $\h$ `energy unstable' at period $X$.
If instead $\E(\h+ \varepsilon v) > \E(\h)$ for all sufficiently small
$\varepsilon > 0$, for each $X$-periodic zero--mean perturbation $v$,
then we call the steady state `energy stable' at period $X$.  (Some
authors call this {\em formal} stability \cite{Holm}.) 
It is conceivable that a steady state might be energy stable and yet
not be a local minimum of the energy.

Our main stability results, in Section~\ref{nonlinear_results}, are
roughly stated as follows.
\begin{itemize}
\item Theorem~\ref{nonlinear_unstable_general}. 
For positive periodic steady states, linear instability
implies energy instability. Hence our linear instability results in
\cite{LP2} imply that every positive $X^*$-periodic steady state is energy
unstable at periods $X=2X^*, 3X^*, \ldots$, and is also energy unstable at period $X^*$ if $g/f$ is a strictly convex function.
\item
Theorems~\ref{nonlinear_unstable_power}--\ref{nonlinear_stable_power}.
Further, for the `power law' coefficients $f(y)=y^n$ and $g(y) = \B
y^m$ with $\B > 0$, we completely characterize energy stability
at period $X^*$ even when $g/f$ is not convex, that is, when $m-n \in (0,1)$.
\end{itemize}
In Section~\ref{relation} we explain how these results for periodic
boundary conditions cover the case of Neumann (`no flux') boundary
conditions as well.

\begin{figure}
  \begin{center}
    \includegraphics[width=.4\textwidth]{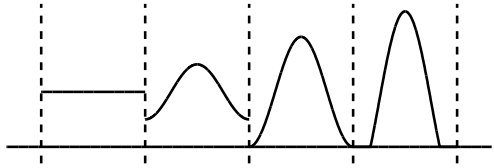}
  \end{center}
  \vspace{-.4cm}  
  \caption{\label{steady} Four types of steady state.}
\vspace{-0.4cm}
\end{figure}

\vskip 6pt 
Then in Section~\ref{energy_steady} we determine the relative energy
levels of three different kinds of steady state: constant steady
states, positive periodic steady states, and zero contact angle
droplet steady states. Figure~\ref{steady} illustrates these three
steady states, as well as showing a nonzero contact angle droplet
steady state (about which we say little in this paper).

We have found it too difficult to determine energy levels when working
with arbitrary coefficients $f$ and $g$, but have obtained fairly
complete answers for power law coefficients.  This provides at least
some insight into the general case. Further insight comes from the
work of Grinfeld and Novick--Cohen \cite{GNC99} on the energy levels
of steady states for the Cahn--Hilliard equation, which has
non--power--law coefficients. The earlier work of Mischaikow \cite{Mis95} applies to a variety of gradient-like bistable equations.

Our main energy level results, in Section~\ref{energy_steady}, for
$f(y)=y^n$ and $g(y) = \B y^m$ are:
\begin{itemize}
\item Theorem~\ref{periodic_constant}. If $m<n$ or $m \geq n+1$ then a
positive periodic steady state always has higher energy than the
constant steady state with the same mean value. For $n < m \leq n +
0.75$ our analytical and numerical work suggest the positive periodic
steady state has lower energy than the constant steady state.
\item Theorem~\ref{2_2_tango}. When $m \approx n + 0.77$ there can be
two steady states with the same period and area, $h_{{\rm ss}1}$ and
$h_{{\rm ss}2}$, with $\min_x h_{{\rm ss}1}(x) < \min_x h_{{\rm
ss}2}(x)$; we essentially prove $h_{{\rm ss}2}$ is energy
unstable and has higher energy than $h_{{\rm ss}1}$, which is 
energy stable.
\item Theorem~\ref{periodic_droplet}. If $m-n \in (-2,0) \cup [1,2)$ then a
positive periodic steady state always has higher energy than a zero
contact angle droplet steady state with the same mean value. 
\item Theorem~\ref{constant_droplet}. The
constant steady state can have higher energy than the zero contact angle droplet steady state with the same mean value. When $m-n \in [1,2)$, for example, a mountain pass scenario can occur, in which the constant steady state is a local minimum of the energy, the positive periodic steady state is an energy unstable saddle, and there is a zero contact angle droplet having lower energy than either of them.
\end{itemize}
For example, Theorem~\ref{periodic_constant} with $m=-1<n=3$ covers the `van der Waals' case 
\begin{equation} \label{Waal}
h_t = - (h^3 h_{xxx})_x - A (h^{-1} h_x)_x 
\end{equation}
with $A>0$. This equation has been studied by a number of other
authors, {\it e.g.} \cite{BBD,WD82,WB99,ZL99}, mostly with regard to
similarity solutions and film rupture (where the solution goes to zero
in finite time, at some point). Our numerical work \cite{LP4} on this equation
pays particular attention to behavior near steady states, and how this
evolves into rupture.


Our energy level results suggest possible basins of attraction around
the stable steady states, and possible heteroclinic connections
between steady states. In the companion article \cite[\S4]{LP4}, we
investigate such possibilities with numerical simulations.  For
example, when $m \approx n + 0.77$ as in Theorem~\ref{2_2_tango}, we
find robust heteroclinic connections between the unstable positive
periodic steady state $h_{{\rm ss}2}$ and the stable one $h_{{\rm
ss}1}$. For the mountain pass scenario in
Theorem~\ref{constant_droplet}, we find that perturbing the saddle
point (the periodic steady state) in one direction leads to relaxation
to the constant steady state and perturbing in the opposite direction
gives apparent relaxation to a droplet. A similar dichotomy was found
for axisymmetric surface diffusion by Bernoff, Bertozzi and Witelski,
\cite[p.\ 744]{BBW}, with perturbed unduloids either relaxing to a
cylinder or else pinching off in finite time.

We also present in \cite[\S5]{LP4} simulations suggesting that small
changes in the `mobility' coefficient $f$ do not break heteroclinic
orbits, but can affect whether or not the solution remains positive as
it evolves.

We discuss some of these conclusions and future directions further in
Section~\ref{conclusions}.

\subsection*{Terminology.} 
We write $\T_X$ for a circle of circumference $X>0$. As usual, one
identifies functions on $\T_X$ with functions on $\R$ that are
$X$-periodic and calls them {\it even} or {\it odd} according to
whether they are even or odd on $\R$.

A positive periodic steady state is assumed to satisfy the steady
state equation classically.  A droplet steady state $\h(x)$ (see
Figures~\ref{steady}c and d) is by definition positive on some
interval $(a,b)$ and zero elsewhere, with $\h \in C^1[a,b]$; we
require $\h$ to satisfy the steady state equation on the open interval
$(a,b)$ only, and to have equal acute contact angles: $0 \leq \h^\d(a)
= - \h^\d(b) < \infty$. (Throughout the paper, if a function has only
one independent variable then we use ${}^\d$ to denote differentiation
with respect to that variable: $\h^\d = (\h)_x$.)

We say a droplet steady state $\h$ has `zero contact angle' if
$0=\h^\d(a)=-\h^\d(b)$, and `nonzero contact angle' otherwise.  A
`configuration' of droplet steady states is a collection of
steady droplets whose supports are disjoint. For more on the steady states
and their properties, see \cite{LP1}.

\section{Energy stability for periodic steady states}
\label{nonlinear_results}

We assume throughout this section that $f(y)$ and $g(y)$ are
$C^2$-smooth for $y>0$, and that $f>0$. Define
\[
r=\frac{g}{f} .
\]
Take $X>0$.

We investigate stability and dynamical questions by means of a
Liapunov energy. A few of the theorems follow directly from our linear
stability results in \cite{LP2}, but most are quite different and
complementary.

\subsection{Definition of the energy, and of energy instability}
\label{instability_defn}
The energy function for
the evolution equation (\ref{evolve}) is defined for $\ell \in
H^1(\T_X)$ to be
\begin{equation} \label{energy_is}
\E(\ell) = \int_0^X \left[ \frac{1}{2} (\ell^\d)^2 - H(\ell) \right] dx ,
\end{equation}
where $H(y)$ is a function with $H^\dd = r = g/f$.

To verify the energy $\E$ is a Liapunov function for the evolution
(\ref{evolve}), suppose $h(x,t)$ is a positive smooth solution
of (\ref{evolve}) that is $X$-periodic in $x$. Bertozzi and Pugh \cite[\S2]{BPLW} observed
(generalizing \cite{oronPRL92,GPS97,oron92}) that $\E$
is dissipated by the evolution: 
\begin{equation*} 
\full{\;}{t} \E(h(\cdot,t)) = - \int_0^X \frac{1}{f(h(x,t))} 
\left[ f(h(x,t)) h_{xxx}(x,t) + g(h(x,t)) h_x(x,t) \right]^2 dx \leq 0. 
\end{equation*}
The dissipation is strict at each time $t$ unless $f(h) h_{xxx} + g(h)
h_x = 0$ for all $x$.  For smooth positive periodic solutions this
occurs only when $h(\cdot,t)$ is a steady state.

Let $\h \in C^4(\T_X)$ be a positive periodic steady state of (\ref{evolve}). It is easy to see (cf.~formula (\ref{var1})) that $\h$ is a critical point for the energy $\E$, with respect to zero-mean perturbations.
\begin{definition}
{\rm Call $\h$ an {\em energy unstable} critical point (with respect to
zero--mean perturbations at period $X$) if there exists a smooth
$X$-periodic perturbation $u(x)$ with mean value zero such that
\[
\E(\h+\varepsilon u) < \E(\h) 
                \qquad \text{for all small $\varepsilon >0$.}
\]
That is, small perturbations in the direction $u$ decrease the
energy. (Some authors call this {\em formal} instability \cite{Holm}.)
}
\end{definition}

An energy unstable steady state is necessarily a saddle point in the energy landscape, since $\E$ is 
{\em increased} by the perturbation $u(x) = \e \cos(2 \pi k x/X)$ for $k \gg 1$.

Energy unstable steady states are not asymptotically stable in
$H^1(\T_X)$, in the following sense: suppose $h(x,t)$ is a positive
smooth solution with initial data $\h+\varepsilon u$; then $h(\cdot,t)
\not \rightarrow \h(\cdot)$ in $H^1(\T_X)$ since for all $t$,
$\E(h(\cdot,t)) \leq \E(h(\cdot,0)) = \E(\h+\varepsilon u) < \E(\h)$
(convergence in $H^1$ would imply convergence in $L^\infty$ and hence
convergence of the energy). In fact, $h(\cdot,t)$ cannot converge to any translate of $\h$, for the same reason.

The last paragraph extends to nonnegative {\it weak} solutions if they
also dissipate the energy.  Weak solutions must sometimes be
considered because solutions of (\ref{evolve}) that are initially
positive might not always remain so, and where they go to zero they
can lose regularity.  See Bertozzi and Pugh~\cite{BPLW,BPFS} for
existence of nonnegative weak solutions that dissipate the energy.

\subsection{Energy instability results} \label{nonlinear_instability}
In \cite[\S2]{LP2} we linearized the evolution equation (\ref{evolve})
around the positive periodic steady state $\h$ and then reduced the
linear stability question to determining the sign of the first
eigenvalue of a certain self-adjoint fourth order linear operator.  We
will not repeat the linearization here, or re-state the linear
stability results of \cite{LP2}.  However we warn readers that when
we say a steady state is `linearly stable', we are including the
neutrally stable case in which the first eigenvalue of the linearized
operator is zero. This is unavoidable: the operator {\em always} has a
zero eigenvalue in its spectrum, corresponding to an infinitesimal
translation of the steady state in space (the evolution equation is translation invariant).

%
%

The next theorem states that if a steady state is linearly unstable
then it is energy unstable.  Also, we present some 
unstable directions, when $g/f$ is strongly convex.
\begin{theorem} \label{nonlinear_unstable_general}
Let $f,g \in C^2(0,\infty)$ with $f>0$.  Take $X>0$ and suppose $\h
\in C^4(\R)$ is an $X$-periodic positive steady state of
(\ref{evolve}).

If $\h$ is linearly unstable with respect to zero--mean perturbations
at period $X$, then it is also energy unstable at period $X$.

In particular, $\h$ is energy unstable at period $X$ if it is
non-constant and either: the least period of $\h$ is $X/j$ for some
integer $j \geq 2$ or else $r=g/f$ is convex ($r^\dd \geq 0$) and
non-constant on the range of $\h$. For example, if $\h$ is
non-constant and $r^\dd>0$, then $\h$ is energy unstable in the
directions $u=\pm \h^\d$ and $\pm \h^\dd$.
\end{theorem}

This is proved in Section \ref{proof_nonlinear_unstable_general}. Note
for example that the theorem covers the van der Waals evolution
(\ref{Waal}), since there $r(y)=A y^{-1}/y^3=A y^{-4}$ is strongly
convex. Thus all positive periodic steady states of the van der Waals
evolution are energy unstable.

These energy unstable steady states (which we observed
above are not asymptotically stable) are presumably {\em nonlinearly}
unstable, in general, but we cannot prove this.
%

\vskip 6pt In Theorem~\ref{nonlinear_unstable_general} we have assumed
$f,g \in C^2(0,\infty)$, which is generally the case for the thin film
equations that are our main motivation. But our arguments are all
local (involving only small perturbations of the steady state), and so
the theorem still holds if the coefficient functions $f(y)$ and $g(y)$
are defined and $C^2$ merely for $y$-values in a neighborhood of 
$[{\h}_{min},{\h}_{max}]$.
\subsection{Review of power law steady states and their rescalings}
\label{power_law_review}
We now turn to {\em power law} coefficients: $f(y) = y^n$
and $g(y)=\B y^m$ for some exponents $n,m \in \R$ and some positive
constant $\B>0$.  Here
\[
r(y) = \B y^{q-1} 
\]
where
\[
\fbox{$q := m-n+1$}
\]
This exponent $q$ determines many properties of the
steady states, 
including (usually) their linear stability.

The evolution equation (\ref{evolve}) becomes
\begin{equation*} 
  h_t = - (h^n h_{xxx})_x - \B(h^m h_x)_x .
\end{equation*}
To state our results on energy stability for this power law evolution,
we first review some properties of the steady states and explain how
to rescale them to solutions of a canonical nonlinear oscillator ODE,
given as equation (\ref{ss4}--\ref{ss5}) below.

\vskip 6pt

We start with a non-constant positive periodic steady state $\h \in
C^4(\T_X)$ of the general evolution (\ref{evolve}). The steady state
condition for (\ref{evolve}) integrates to give $f(\h) \h^\ddd + g(\h)
\h^\d = C$ for some constant $C$. The least period of $\h$ is $X/j$
for some integer $j \geq 1$.

One finds that the constant $C$ (the flux) equals zero, by dividing
$f(\h) \h^\ddd + g(\h) \h^\d = C$ by $f(\h)>0$ and integrating over a
period (cf.~\cite{oron92,oron94}).  Hence the steady state satisfies
\begin{equation} \label{ss1}
\h^\ddd + r(\h) \h^\d = 0.
\end{equation}
[If $\h$ were a droplet steady state then again $C = 0$, by
\cite[Theorem~2.1]{LP1}, and equation (\ref{ss1}) would hold wherever
$\h$ is positive.]

Integrating, the steady states have a nonlinear oscillator
formulation:
\begin{equation} \label{ss2}
  \h^\dd + H^\d(\h) = 0
\end{equation}
holds wherever $\h$ is positive.  Here $H(y)$ is a function with
$H^\dd=r=g/f$; if we regard $x$ as a `time' variable then $\frac{1}{2} \h^\d(x)^2 + H(\h(x))$ is a conserved quantity.

Returning to the power law evolution, remember $r(y)=\B y^{q-1}$. Thus
for $q \neq 0$ we can write the steady state equation (\ref{ss2}) as
\begin{equation} \label{ss3}
  \h^\dd + \frac{\B \h^q - D}{q} = 0
\end{equation}
for some constant $D$.
For $q=0$ the analogous equation is $\h^\dd + \B \log \h - D = 0$.  
This oscillator equation involves three constants: $q$, $\B$, and
$D$. We remove $\B$ and $D$
by rescaling: let
\renewcommand{\arraystretch}{2}
\begin{equation}  \label{rescaling}
k(x) = 
\left\{
        \begin{array}{rl}
        \left( \frac{\B}{D} \right)^{\! 1/q} 
        \h \! \left( \left( \frac{D}{\B} \right)^{\! 1/2q}
        \frac{x}{D^{1/2}} \right) , & q \neq 0, \\
        e^{-D/\B} \h \! \left( e^{D/2\B} \frac{x}{\B^{1/2}} \right) , 
        & q = 0.
        \end{array}
\right.
\end{equation}
\renewcommand{\arraystretch}{1}
\noindent 
For $q \neq 0$ this rescaling uses that $D > 0$, by
\cite[\S3.1]{LP1}. [A different rescaling would be used to study
droplet steady states \cite[\S\S3.2,4]{LP1}.] The steady state
equation (\ref{ss3}) rescales to
\begin{eqnarray} \label{ss4}
k^\dd + \frac{k^q - 1}{q} & = & 0, \qquad q \neq 0, \\ 
\label{ss5}
        k^\dd + \log k            & = & 0, \qquad q   =  0 .
\end{eqnarray}
Differentiating, we find for all $q$ that $k^\ddd + k^{q-1} k^\d = 0$, 
and so $k$ satisfies $\left( k^n k^\ddd + k^m k^\d \right)^\d = 0$, {\it i.e.} it is a steady state of $k_t = - (k^n k_{xxx})_x - (k^m k_x)_x$.

Since $\h$ is non-constant, positive and 
periodic, we know $k^\dd(x_0) >
0$ for some point $x_0$. Evaluating 
(\ref{ss4}--\ref{ss5}) at $x_0$ shows the minimum value of $k$ is less than $1$. Also $k^\d(0)=0$ since (after a suitable translation) $\h$ has its minimum at $x=0$. 
Introducing the notation $k_\alpha$ for the solution $k$ that has
minimum value $\alpha \in (0,1)$, at $x=0$, we have
\begin{equation} \label{init_cond}
0 < k_{\alpha}(0)=\alpha < 1, \qquad k_\alpha^\d(0) = 0 .
\end{equation}

\begin{figure}
  \begin{center}
                \includegraphics[width=.3\textwidth]{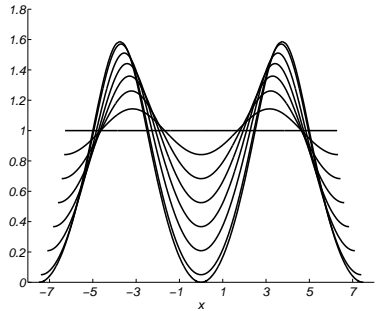}
  \end{center}
  \vspace{-.2cm}  
  \caption{\label{ssq3} 
        Steady states $k_\alpha(x)$, when $q = 3$.}
\end{figure}

Thus every steady state $\h$ can be rescaled to a $k_\alpha$, as
above.  Conversely, for each $q \in \R$ and $\alpha \in (0,1)$ there
exists a unique smooth positive periodic $k_\alpha$ satisfying
equations (\ref{ss4}--\ref{ss5}) and (\ref{init_cond}) (see
\cite[Proposition~3.1]{LP1}). The same holds
for $\alpha=0$ when $q>-1$, except that $k_0$ may be only $C^1$-smooth
at $x=0$ (see \cite[Theorem~3.2]{LP1}).  To illustrate,
Figure~\ref{ssq3} plots the steady states $k_\alpha$ over two periods,
for $q=3$ and eight $\alpha$-values between $0$ and $1$; see
\cite[\S6.1]{LP1} for details.

Note that the map $(\alpha,x) \mapsto k_\alpha(x)$ is
$C^\infty$-smooth for $(\alpha,x) \in (0,1) \times \R$, by an ODE
theorem giving smooth dependence on the initial data \cite[Ch.~V
\S4]{Hartman}.  We write
\[
        P=P(\alpha) \qquad \text{and} \qquad A=A(\alpha)
\]
for the least period of $k_\alpha$ and for the area under its graph, $A=\int_0^P k_\alpha(x) \, dx$, respectively. Then $P$ and $A$ are smooth
functions of $\alpha$ that approach $2\pi$ as $\alpha \rightarrow 1$, by \cite[Lemma~6]{LP2}. The function
\[
E(\alpha) := P(\alpha)^{3-q} A(\alpha)^{q-1} = P(\alpha)^2 [A(\alpha)/P(\alpha)]^{q-1}
\]
determines whether the steady state is energy unstable or stable, in several results below.  

\vskip 6pt 

The above rescaling ideas are a useful tool throughout the paper. We
hope this tool does not obscure the fact that stability and energy
level properties for equations of type (\ref{evolve}) seem to be
determined by the period map of a family of steady states $\h$ with
fixed area but varying amplitudes and periods, or alternatively the
area map of a family of steady states with fixed period but varying
amplitudes and areas. This is how one should think of the function
$E(\alpha)$; see \cite[\S6.3]{LP2} for more on this. The same
underlying idea appears in the work of Grinfeld and Novick--Cohen
\cite{GNC99} on the Cahn--Hilliard equation, an equation which is not
amenable to rescaling in the same way.

\subsection{Energy in/stability for the power law evolution}
\label{power_stability}
\begin{theorem} \label{nonlinear_unstable_power}
Let $\h \in C^4(\R)$ be a non-constant positive $X$-periodic steady
state of the power law evolution $h_t = -(h^n h_{xxx})_x - \B (h^m
h_x)_x$.
Translate $\h$ to put its minimum at $x=0$ so that $\h$ rescales to
$k_\alpha$ for some $\alpha \in (0,1)$, as in \S \ref{power_law_review}.
 
\noindent If $q < 1$ or $q > 2$ then $\h$ is energy unstable
in the directions $u=\pm \h^\d$ and $u=\pm \h^\dd$.

\noindent If $q = 2$, or if $q > 1$ and $E^\d(\alpha) > 0$, then $\h$ is energy unstable.
\end{theorem}
We prove the theorem in Section
\ref{proof_nonlinear_unstable_power}. Its first statement follows
immediately from Theorem~\ref{nonlinear_unstable_general}, since
$r(y)=\B y^{q-1}$ is strongly convex ($r^\dd>0$) when $q<1$ or $q >
2$. The final statement of the theorem certainly applies when $q \geq
2$, since then $E^\d>0$ by \cite[Theorem~11]{LP2}.  Computational
studies \cite[\S6.1]{LP1} suggest $E^\d(\alpha) > 0$ for all $\alpha
\in (0,1)$ except when $q \in [1,1.794]$ (approximately);
Figure~\ref{E_close-up} plots $E(\alpha)$ for certain of these
$q$-values.

To help explain the appearance of the criterion $E^\d(\alpha) > 0$ for
energy instability, in the theorem, we refer the reader to
\cite[\S6.3]{LP2} and to \cite[\S3]{LP4} for a `bifurcation diagram'
interpretation of the function $\alpha \mapsto E(\alpha)$ in terms of
a family of steady states $\h$ with fixed area but varying amplitudes
and periods, or fixed period but varying amplitudes and areas.

\vskip 6pt The case $1<q<2$ with $E^\d(\alpha) \leq 0$ is not covered
by Theorem~\ref{nonlinear_unstable_power}. By \cite[Theorem~9]{LP2} we
know $\h$ is linearly stable in this case. We have been unable to
prove {\em non}linear stability, but we do prove in
Theorem~\ref{nonlinear_stable_power} below that if $E^\d(\alpha) < 0$
then small perturbations of $\h$ in every possible direction strictly
increase the energy; this is consistent with nonlinear stability.
\begin{definition}
  {\rm Let $\h \in C^4(\T_X)$ be a non-constant positive periodic
  steady state of (\ref{evolve}). Call $\h$ {\em energy stable} (with
  respect to zero--mean perturbations at period $X$) if for each $u
  \in H^1(\T_X) \setminus \{ 0 \}$ with $\int_0^X u \, dx = 0$ we have
\[
\E(\h+\varepsilon u) > \E(\h) 
                \qquad \text{for all small $\varepsilon >0$.}
\]
\noindent (Some authors would call this {\em formal} stability \cite{Holm}.)
}
\end{definition}

A steady state might perhaps be energy stable without being a local
minimum of the energy.  (For example, the function $f(x,y) = (y-x^2)
(y - 3 x^2)$ on $\R^2$ has the property that origin is a local minimum
point on each straight line through the origin, though it is not a
local minimum in $\R^2$.) Another cautionary note is that the energy
is insensitive to translation: in particular, a steady state and its
translates have the same energy, consistent with the translation
invariance of the evolution equation itself.  Numerical simulations in
our article \cite[\S4.4]{LP4} demonstrate that perturbations of an
energy stable steady state can evolve towards a translate of that
steady state.  This suggests that any asymptotic or nonlinear
stability result that can be proved will hold only up to translation.

In Section \ref{proof_nonlinear_stable_power} we prove:
\begin{theorem} \label{nonlinear_stable_power}
Let $\h \in C^4(\R)$ be a non-constant positive periodic steady state
of the power law equation $h_t = -(h^n h_{xxx})_x - \B (h^m h_x)_x$,
with $\h$ having least period $X$. Translate $\h$ to put its minimum
at $x=0$ so that $\h$ rescales to $k_\alpha$ for some $\alpha \in
(0,1)$, as in \S \ref{power_law_review}.
 
If $1 < q < 2$ and $E^\d(\alpha) < 0$ then $\h$ is energy stable.
\end{theorem}
The hypothesis $E^\d(\alpha) < 0$ seems numerically to hold 
for all $\alpha \in (0,1)$
when $1 < q \leq 1.75$, as indicated by Figure~3 in \cite{LP2}.

We have no energy stability result when $1<q<2$ and $E^\d(\alpha) = 0$; fortunately it seems $E^\d(\alpha) = 0$ for at most {\em one} $\alpha$-value, for each $q$, as shown numerically by Figures~3--5 in \cite{LP2}. Theorem~9 of
\cite{LP2} does imply linear stability when $E^\d(\alpha) = 0$, but
\cite[Theorem~10(b)]{LP2} shows the space of neutrally stable
directions is two dimensional (rather than one dimensional as when
$E'(\alpha) < 0$) and this might perhaps lead to instability.

We next address the $q=1$ case.
\begin{lemma} \label{q_one_nonlinear}
Let $q=1$ ({\it i.e.} $m=n$) and suppose $\h \in C^4(\R)$ is a
non-constant positive periodic steady state of $h_t = - (h^n
h_{xxx})_x - \B (h^n h_x)_x$ with least period $X$, and translate so
that $\h$ has its minimum at $x=0$. Then $\h$ is not asymptotically
stable in $H^1(\T_X)$ with respect to even perturbations.
\end{lemma}
\begin{proof}
The steady state equation (\ref{ss3}) with $q=1$ has general solution
$\h(x)=D/\B + c \cos (\sqrt{\B}x)$, where we have used that $\h$ has
an extremum at $x=0$. Hence the period is $X=2\pi/\sqrt{\B}$, and for
all small $\varepsilon$ we see that the perturbed function $\h(x) +
\varepsilon \cos (\sqrt{\B}x)$ is another positive periodic steady
state solution. Thus $\h$ is not asymptotically stable in $H^1(\T_X)$
with respect to even perturbations.
\end{proof}

The steady state is of course not asymptotically stable with respect
to {\em general} perturbations, since one can always perturb by
translating the steady state a small distance (this remark applies for
all $q$). However we are really not interested in such simple
translational motion.
Also,
translational perturbations are not even permitted under the Neumann
boundary conditions that we consider below, in Section~\ref{relation}.

To summarize, when $q=1$ the positive periodic steady states are
linearly neutrally stable with respect to zero--mean perturbations of
the same period, by \cite{GPS97} or \cite[Lemma~8]{LP2}, and are not
asymptotically stable with respect to `even' perturbations, by the
above lemma. Our numerical simulations in \cite[\S4.3]{LP4} and those
of \cite{GPS97} suggest that a wide range of small perturbations yield
solutions relaxing to nearby positive periodic steady states, suggesting
they are nonlinearly stable.

\vskip 6pt 

In the companion paper \cite[\S4]{LP4} we illustrate the above
stability and instability theorems for the power law evolution with a
variety of numerical simulations.  There we find not only the short
time behavior suggested by the energy (in)stability results, but also
some longer time limits that are suggested by the energy level results
in Section~\ref{energy_steady}.

\subsection{Odd perturbations}
\label{odd_global}
Returning momentarily to general coefficients $f$ and $g$, in Section
\ref{proof_nonlinear_odd_power} we prove the energy increases under
{\em odd} perturbations, when $r$ is concave.
\begin{theorem} \label{nonlinear_odd_power}
Let $\h \in C^4(\R)$ be a non-constant positive periodic steady state
of (\ref{evolve}) with coefficient functions $f,g \in C^2(0,\infty),
f>0$. Suppose $\h$ has least period $X$, and translate $\h$ to put its
minimum at $x=0$.

If $r=g/f$ is strongly concave ($r^\dd < 0$) then for every  
nontrivial $u \in H^1(\T_X)$
that is odd and is such that $\h+u>0$, we have $\E(\h+u) > \E(\h)$.
\end{theorem}
The theorem is global since the perturbations are not required to be
small, and it is consistent with asymptotic stability (although
convergence to a {\em translate} of $\h$ seems more likely than
convergence to $\h$ itself).

Theorem~\ref{nonlinear_odd_power} applies in the power law case with
$1 < q < 2, r(y)=y^{q-1}$.  Another example with strongly concave $r$
is an equation \cite[eqn. (24)]{Lewis94} describing the dynamics of a
population of aphids, for which $f(y)=y, g(y)=y-c$ and $r(y)=1-c/y$.

\subsection{Relation between the periodic and Neumann stability problems}
\label{relation}

Suppose $\h$ is an even $X$-periodic steady state of the evolution
equation (\ref{evolve}) with extrema at $x=0,\pm X/2,\ldots$, so that
$\h^\d=\h^\ddd=0$ at these points. As we observed at the end of
\cite[\S2.5]{LP2}, linear instability of $\h$ with respect to periodic
boundary conditions on $(-X/2,X/2)$ is equivalent to linear
instability with respect to Neumann (`no flux') conditions on the
half-interval $(0,X/2)$; these Neumann conditions are: $h_x=h_{xxx}=0$
at $x=0,X/2$.

The energy of a positive smooth solution is still dissipated in the
case of Neumann boundary conditions, and obviously energy instability
of $\h$ in an {\em even} direction $u(x)$ for the periodic problem on
$(-X/2,X/2)$ is equivalent to energy instability in the direction
$u(x)$ for the Neumann problem on $(0,X/2)$.  (If the perturbation
$u(x)$ is even and has mean value zero on $(-X/2,X/2)$ then it has
mean value zero on $(0,X/2)$ as well.) Similarly, `periodic' energy
stability in all even directions on $(-X/2,X/2)$ is equivalent to
`Neumann' energy stability in all directions on $(0,X/2)$.

In short, for the Neumann problem on $(0,X/2)$, the stability result
in Theorem~\ref{nonlinear_stable_power} still holds, and the
instability claims involving $\pm \h^\dd$ in
Theorems~\ref{nonlinear_unstable_general} and
\ref{nonlinear_unstable_power} also still hold, since these are even
functions; the claims involving $\pm \h^\d$ do not carry over, since
those are odd.

\section{Relative energy levels of periodic, constant and droplet steady states}
\label{energy_steady}

In this section we investigate the phase space of the power law
equation $h_t = -(h^n h_{xxx})_x - \B (h^m h_x)_x$ by comparing the
value of the energy at the positive periodic, constant and zero-angle
droplet steady states. Let $X>0$ and recall $q=m-n+1$.

\subsection{Positive periodic vs.~constant steady states}
\label{energy_positive}
The fluid volume $\int_0^X h(x,t) \, dx$ is conserved by the
evolution, under periodic boundary conditions, and so the mean value
\[
\overline{h} := \frac{1}{X} \int_0^X h \, dx
\] 
is constant in time.  Suppose the initial data $h(\cdot , 0)$ arises
from a small zero-mean perturbation of $\h$.  It is natural to ask
whether $h$ can converge (while staying positive and smooth) towards
the constant steady state $\overline{\h}$, as $t \rightarrow
\infty$. This {\em cannot} happen if $\E(\h) < \E(\overline{\h})$ and
also $h(\cdot,0)$ is close enough to $\h$ so that $\E(h(\cdot,0)) <
\E(\overline{\h})$, because the energy is dissipated by the evolution.
\begin{theorem} \label{periodic_constant}
Let $\h \in C^4(\R)$ be a non-constant positive periodic steady state
of the power law equation $h_t = -(h^n h_{xxx})_x - \B (h^m h_x)_x$,
with least period $X$. Translate $\h$ to put its minimum at $x=0$ so
that $\h$ rescales to $k_\alpha$ for some $\alpha \in (0,1)$, as in
\S\ref{power_law_review}.

If $q \geq 2$ or $q < 1$ then $\E(\h) > \E(\overline{\h})$.

If $q=1$ then $\E(\h) = \E(\overline{\h})$.

If $1<q<2$ and $E^\d>0$ on $(\alpha,1)$ then $\E(\h) > E(\overline{\h})$.

If $1<q<2$ and $E^\d<0$ on $(\alpha,1)$ then $\E(\h) < E(\overline{\h})$.
\end{theorem}
The theorem is proved in Section~\ref{proof_periodic_constant}. For $q=-3$ it was observed numerically in \cite[\S3]{WB00} that $\E(\h) > \E(\overline{\h})$.

When $1 < q \leq 1.75$, numerical evidence in Figure~3 of \cite{LP2}
suggests $E^\d(\alpha)<0$ for all $\alpha$. If this is true, then
$\E(\h) < \E(\overline{\h})$ by the theorem and so there can be no heteroclinic connection from $\h$ to the constant steady state $\overline{\h}$; indeed in \cite[\S4.4]{LP4} we find numerically for $q=1.5$ that $\h$ is asymptotically stable, up to translation.  

For $1.795 \leq q <2$,
Figure~5 of \cite{LP2} suggests $E^\d(\alpha)>0$ for all $\alpha$ and so
$\E(\h) > \E(\overline{\h})$ by the theorem.
Thus when $q \geq 1.795$ or $q<1$, the instability result Theorem~\ref{nonlinear_unstable_power} and the energy level result Theorem~\ref{periodic_constant} together lead us to suspect the existence of a heteroclinic connection from $\h$ to $\overline{\h}$.
\begin{figure}[h] 
   \begin{minipage}{.5\linewidth} 
    \begin{center}
    \includegraphics[width=.7\textwidth]{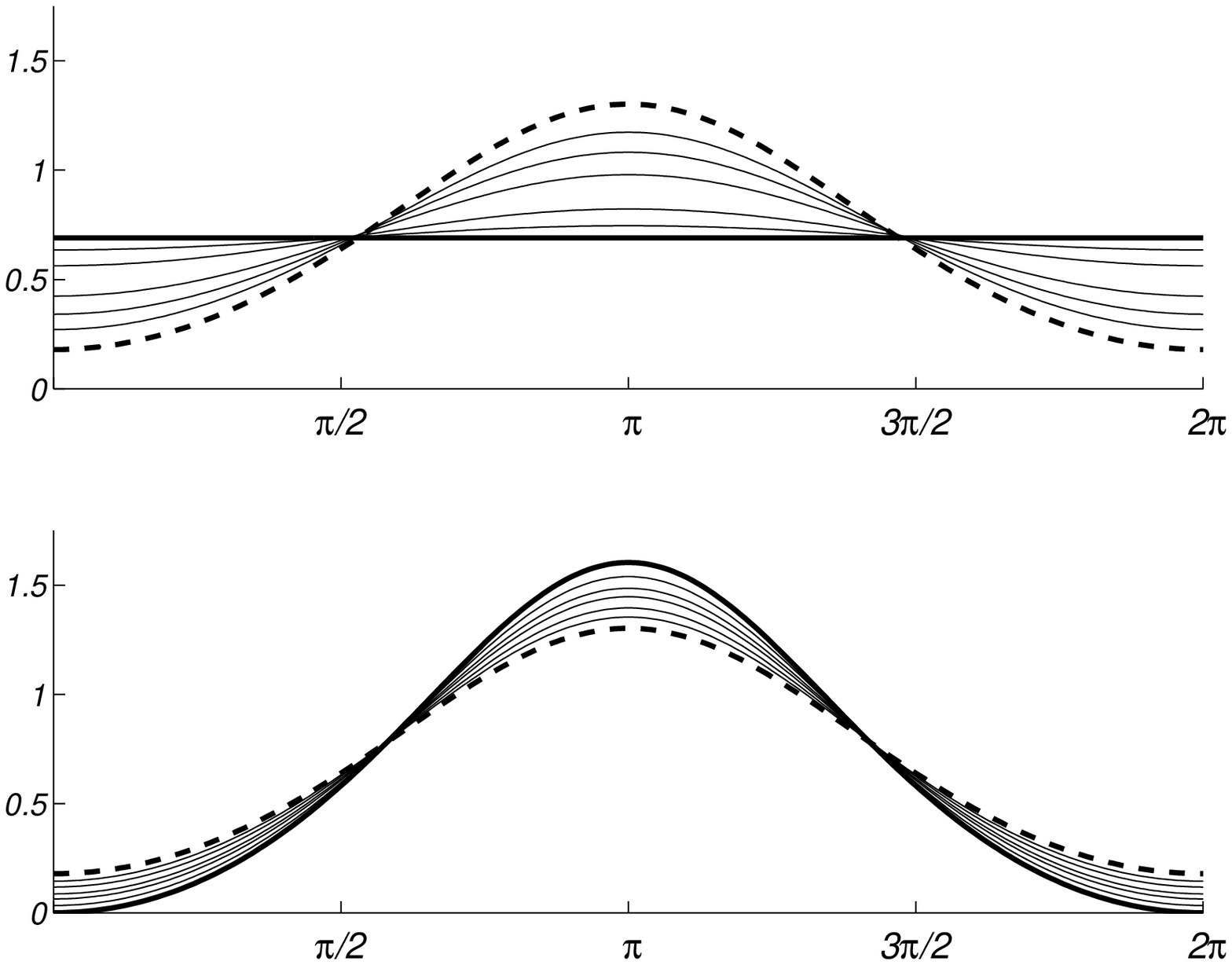}
    \vspace{-.2cm}      
    \caption{\label{q2.5_fig2}$q = 2.5$ and $n=1$, dashed line: initial data.}
    \end{center}
  \end{minipage}%
 \begin{minipage}{.5\linewidth} 
    \begin{center}      
    \includegraphics[width=.7\textwidth]{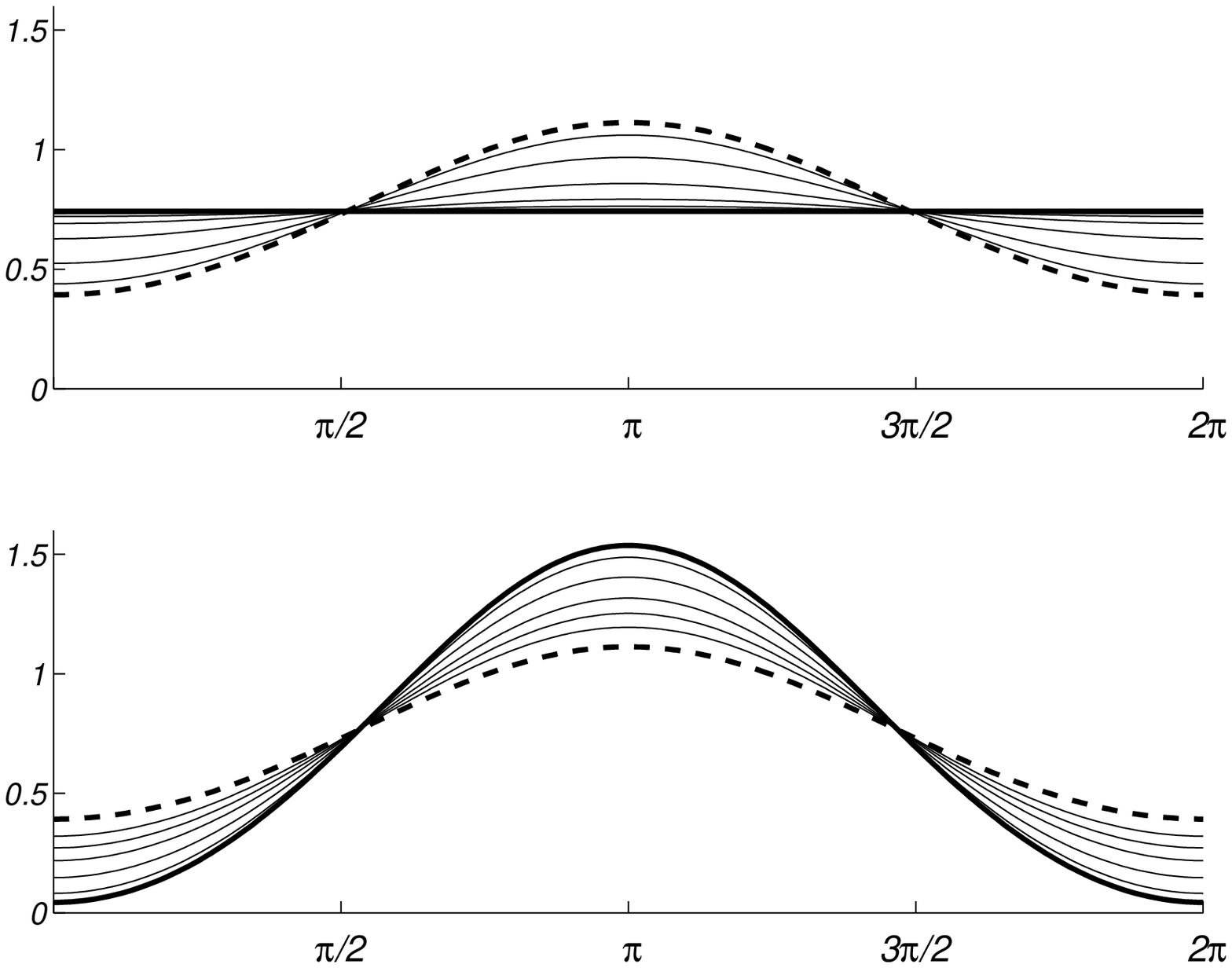}
    \vspace{-.2cm}      
    \caption{\label{q1.768_het}$q = 1.768$, $n=1$, dashed line: initial data.}
    \end{center}
  \end{minipage}
  \vspace{-0.2cm}
\end{figure}
In Figure~\ref{q2.5_fig2} we present numerical simulations of such
heteroclinic orbits, taken from the companion article
\cite[\S4.6]{LP4}.  The top part of the figure presents an orbit
connecting the positive periodic steady state to the constant steady
state.  The bottom part relates to \S\ref{energy_periodic_droplet} and
presents a solution connecting a perturbation of the same positive
periodic steady state in finite time to a droplet profile (which will
not in general be a steady state).

In particular the van der Waals equation (\ref{Waal}) has
$q=-1-3+1=-3$, and so by the first claim of
Theorem~\ref{periodic_constant}, the energy of any positive periodic
steady state is greater than that of the constant steady state. This
was observed numerically by Witelski and Bernoff \cite[\S3]{WB00}. See
\cite{WD82,WB00} and our companion paper \cite[\S4.1]{LP4} for
numerical simulations of this equation and a discussion of the
literature.

For $q \in (1.75,1.794]$, see the remarks around Theorem~\ref{2_2_tango}.

For the Cahn--Hilliard equation, the analogue of
Theorem~\ref{periodic_constant} (comparing energy levels of
non-constant and constant steady states) can be found in Grinfeld and
Novick--Cohen's work \cite[Theorem 4.1]{GNC99}; further,
\cite[\S7]{GNC99} discusses a number of results on existence of
heteroclinic connections. See also \cite[\S3.4]{Mis95}.

\subsection{Positive periodic vs.~droplet steady states}
\label{energy_periodic_droplet}
We do not yet have stability results for {\em droplet} steady
states. But here we show under certain conditions that the energy of a zero-angle droplet steady state
must be lower than that of a positive periodic steady state $\h$ whose
period exceeds the length of the droplet. In our article \cite[\S4]{LP4} we show numerically that in these cases,
the droplet seems to be strongly attracting, with heteroclinic orbits
from $\h$ towards the droplet steady state as shown in the bottom part
of Figure~\ref{q2.5_fig2}.
\begin{theorem} \label{periodic_droplet}
  Let $\h \in C^4(\R)$ be a non-constant positive periodic steady state 
  of the power law equation $h_t = -(h^n h_{xxx})_x - \B (h^m h_x)_x$,
  with $\h$ having least period $X$ and area $\Ass=\int_0^X \h \, dx$.

  If $-1<q<1$ or $2 \leq q < 3$, or if $1<q<2$ and $E^\d(\alpha)>0$ for all $\alpha \in (0,1)$, then there exists a zero contact angle droplet steady state
  $\hhat$ with length $\hat{X} < X$ and area $\hat{A}=\Ass$.  Furthermore
  $\E(\hhat)<\E(\h)$.
\end{theorem}
We prove the theorem in Section~\ref{proof_periodic_droplet}.

\noindent \textsc{Remarks}

1. Steady states with zero contact angle can occur only for $q>-1$, by
\cite[\S2.2]{LP1}. Hence we do not consider $q \leq -1$ in the
theorem; indeed variational methods show there can be droplet steady
states that do not have acute contact angles \cite{ANC93}.

2. The theorem presumably applies to $1.795 \leq q < 2$, since it
seems numerically that $E^\d > 0$ for those $q$-values.  For $q \geq
3$ we think a similar theorem might hold but with $\hhat$ being a {\em
configuration} of disjoint zero--angle droplets.

3. Since $\E(\h) > \E(\hhat)$, there might be an orbit from $\h$ to
$\hhat$.  This orbit might describe a positive solution that converges
to the nonnegative droplet profile as $t \to \infty$, or it might
describe a positive solution that loses positivity in finite time and
{\it then} approaches the droplet as a nonnegative weak solution.  For
$-1<q<1$ and $2 \leq q < 3$, Theorems~\ref{periodic_constant} and
\ref{periodic_droplet} are consistent with the idea that the unstable
positive periodic steady state $\h$ and its stable manifold form a
separatrix between the basin of attraction of the constant steady
state and the basin of attraction of the droplet steady state.  In
particular, after perturbing $\h$ in one direction one seems to find a
solution that converges to the constant steady state, while perturbing
in the opposite direction often yields a solution that converges to a
droplet profile.  We present some numerical evidence for such behavior
in Figure~\ref{q2.5_fig2}, and discuss this at length in the companion
article \cite[\S4]{LP4}. See also the `mountain pass' remark after
Theorem~\ref{constant_droplet}.

\vskip 6pt 
The preceding theorem and its Remarks are definitely not valid for $q
\in (1,1.75]$, since for these $q$-values we show there does not even
{\em exist} a zero contact angle steady state with length less than
$X$:
\begin{theorem} \label{nonzero_angle}
  Let $1 < q \leq 1.75$ and suppose $E^\d(\alpha) < 0$ for all $\alpha
  \in (0,1)$. Let $\h \in C^4(\T_X)$ be a non-constant positive
  periodic steady state of the power law equation $h_t = -(h^n
  h_{xxx})_x - \B (h^m h_x)_x$, with $\h$ having least period $X$ and
  area $\Ass=\int_0^X \h \, dx$. Assume $\hhat$ is nonnegative and
  piecewise-$C^1$ on $\T_X$, has area $\Ass$, and is smooth on the set
  where it is positive and satisfies there the `nonlinear oscillator'
  steady state equation (\ref{ss3}).

  Then $\hhat$ is either constant or is a translate of $\h$, or is a
  configuration of nonzero contact angle droplet steady states. 
  Specifically,
  $\hhat$ cannot be a zero contact angle droplet steady state.
\end{theorem}
We prove the theorem in Section~\ref{proof_nonzero_angle}.  Note that
the hypothesis $E^\d < 0$
seems to hold for $1<q \leq 1.75$, by the numerical evidence in
Figure~3 of \cite{LP2}.

\begin{figure}
    \begin{center}
    \includegraphics[width=.3\textwidth]{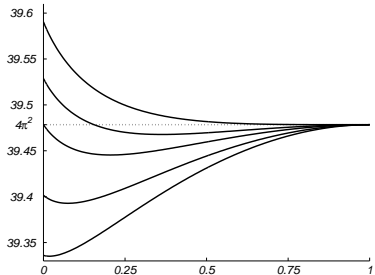}
  \vspace{-.2cm}  
  \caption{\label{E_close-up} 
    Plots of $E(\alpha)$ for 
    $q = 1.75, 1.76, 1.768, 1.78,1.79$.
    Top curve: $q=1.75$, the curves move down as $q$ increases.}
    \end{center}
\end{figure}

\vskip 6pt 
Finally, for $q \in (1.75,1.794]$ approximately, we know by the
analytical and numerical work in \cite[\S5.1]{LP1} that there can be
{\em two} positive periodic steady states with the same period and
area.  The next theorem shows that the steady state with smaller
minimum value (and larger amplitude) is energy stable, and has lower
energy than the other one, which is energy unstable.  This leads us to
suspect there exists an orbit from the energy unstable steady state to
the energy stable one, at least when the steady states have been
chosen to have their minima at the same location.  The bottom plot of
Figure~\ref{q1.768_het} presents a numerical simulation of such an
orbit, taken from \cite[\S4.5]{LP4}. The top plot of the figure
presents an orbit connecting the unstable positive periodic steady
state to the constant steady state, similar to the top plot of
Figure~\ref{q2.5_fig2}.
\begin{theorem} \label{2_2_tango}
  Assume $1<q<2$ and there exists $\alpha_{crit} \in (0,1)$ with
  $E^\d(\alpha) < 0$ on $(0,\alpha_{crit})$ and $E^\d(\alpha) > 0$ on
  $(\alpha_{crit},1)$, and assume $\alpha \mapsto \alpha
  P(\alpha)^{2/(q-1)}$ is strictly increasing for $\alpha \in (0,1)$.

  Suppose $h_{\rm ss1}$ and $h_{\rm ss2}$ are non-constant positive
  periodic steady states of the power law equation $h_t = -(h^n
  h_{xxx})_x - \B (h^m h_x)_x$, with $h_{\rm ss1}$ and $h_{\rm ss2}$
  having the same least period $X$ and same area $\int_0^X h_{\rm ss1}
  \, dx = \int_0^X h_{\rm ss2} \, dx$.

If $h_{\rm ss1}(x)$ has lower minimum value than $h_{\rm ss2}(x)$,
then $h_{\rm ss1}$ is energy stable, $h_{\rm ss2}$ is energy unstable,
and $\E(h_{\rm ss1}) < \E(h_{\rm ss2})$. Furthermore, $\E(h_{\rm ss2})
> \E(\overline{h_{\rm ss2}})$.
\end{theorem}
We prove this in Section~\ref{proof_2_2_tango}. The hypothesis about
$E(\alpha)$ being first strictly decreasing and then strictly
increasing is confirmed numerically for $q$ in the interval
$(1.750,1.794]$ by \cite[\S6.1]{LP1}; see Figure~\ref{E_close-up}.
Numerical work also confirms the hypothesis about $\alpha
P(\alpha)^{2/(q-1)}$ being a strictly increasing function of $\alpha$
for all $q \in (1,2)$.

Theorem~\ref{2_2_tango} is analogous to \cite[Theorem 4.1(v)]{GNC99} for the transitional and metastable cases of the Cahn--Hilliard equation, where Grinfeld and Novick--Cohen show the energy of a monotonic `interface' solution is less than that of a monotonic `spike' solution having the same length and area.

\subsection{Constant steady states: stability results}
\label{energy_constant_stability}
For any number $\hbar > 0$ the constant function $\h \equiv
\hbar$ is a steady state of the general evolution equation
(\ref{evolve}). We now develop analogues for this constant steady
state of our earlier stability results. Then in the following
subsection we can compare the energies of the constant and droplet
steady states, for power law coefficients.

We start by recalling a characterization of linear instability for
steady states. From \cite[Appendix A]{LP2}, an $X$-periodic positive
steady state $\h \in C^4(\R)$ of (\ref{evolve}) is linearly unstable
with respect to zero--mean perturbations at period $X$ if and only if
$\tau_1(\h)<0$, where
\begin{equation} \label{rayleigh}
\tau_1(\h) 
:= \min_u \frac{\int_0^X \left[ (u^\d)^2 - r(\h) u^2 \right] dx}{\int_0^X u^2 \, dx};
\end{equation}
the minimum here is taken over nonzero $u \in H^1(\T_X)$ with
$\int_0^X u \, dx = 0$. The trial function $u$ corresponds to a
perturbation of $\h$.
\begin{theorem} \label{constant_stability}
  Let $f,g \in C^2(0,\infty)$ with $f>0$, and write $r=g/f$. Let
  $\hbar,X>0$. Then for the constant steady state $\h \equiv \hbar$ of
  (\ref{evolve}), the eigenvalue $\tau_1$ in (\ref{rayleigh}) is
  $\tau_1(\hbar) = (2\pi/X)^2 - r(\hbar)$. The
  $\tau_1(\hbar)$-eigenspace is spanned by $\sin(2\pi x/X)$ and
  $\cos(2\pi x/X)$.

  Thus with respect to zero--mean perturbations at period $X$, the
  constant steady state is
\begin{eqnarray*}
\text{linearly unstable} & \text{if} & r(\hbar) X^2 > 4 \pi^2 , \\
\text{linearly neutrally stable} & \text{if} & r(\hbar) X^2 = 4 \pi^2 , \\
\text{linearly asymptotically stable} & \text{if} & r(\hbar) X^2 < 4 \pi^2 .
\end{eqnarray*}

(a) If $r(\hbar) X^2 > 4 \pi^2$, or if $r(\hbar) X^2 = 4 \pi^2$ and $r^\dd(\hbar)>0$, then the constant steady state is energy unstable
in the directions $\pm \sin(2\pi x/X)$ and $\pm \cos(2\pi x/X)$.

(b) If $r(\hbar) X^2 < 4 \pi^2$, or if $r(\hbar) X^2 = 4 \pi^2$ and
$r^\dd(\hbar)<0$, then the constant steady state is energy stable with
respect to zero--mean perturbations of period $X$.  In fact, if
$r(\hbar) X^2 < 4 \pi^2$ then the constant steady state is a strict
local minimum of the energy with respect to zero--mean perturbations
in $H^1(\T_X)$, and $\overline{h}$ is nonlinearly stable under the
evolution (\ref{evolve}), in an $H^1(\T_X)$-sense made precise in the
proof.
\end{theorem}
We prove the theorem in Section~\ref{proof_constant_stability}. Its
linear stability assertions are well-known and
are included for the sake of completeness.

Goldstein, Pesci and Shelley \cite[\S{IIIB}]{GPS97} used the energy to
prove nonlinear instability of the constant steady state for the $q=1$
case ($f(y) = y^n, g(y) = \B y^n, X=2\pi$) with either $2 \leq \B < 4$
or $\B = j^2$ for some integer $j \geq 2$.
\subsection{Constant vs.~droplet steady states}
\label{energy_constant_droplet}

Consider power law coefficients, so that $r(y)=\B y^{q-1}$; then the
previous theorem shows that the stability of the constant steady state
$\hbar$ is determined by whether the quantity $\B \overline{h}^{q-1}
X^2$ is $>,=,< 4 \pi^2$.

Fix $X>0$. Does a zero-angle droplet steady state exist with length at
most $X$ and with the same area $\hbar X$ as the constant steady
state?  If such a droplet steady state exists, can it have lower
energy than the constant steady state?

In this direction, in Section~\ref{proof_constant_droplet} we prove:
\begin{theorem} \label{constant_droplet}
Let $\hbar,X>0$, and consider the constant steady state $\h \equiv
\hbar$ of the power law evolution equation $h_t = -(h^n h_{xxx})_x -
\B (h^m h_x)_x$.

  (a) Suppose $-1 < q < 3$. Then there exists a zero contact angle
  droplet steady state $\hhat$ of length $\hat{X} \leq X$ and area
  $\hbar X$ if and only if
\begin{equation} \label{cd1}
  \B \hbar^{q-1} X^2 \geq E(0) =: E_0(q) .
\end{equation}
If such a droplet steady state exists, then $\E(\overline{h}) > \E(\hhat )$ if and only if
\begin{equation} \label{cd2}
  \B \hbar^{q-1} X^2 > A(0)^2 \left[ \frac{3+q}{(3-q)(q+1)} \right]^{(3-q)/q}
  =: L(q) \qquad \text{(for $-1<q<3,q \neq 0$)}
\end{equation}
or $\B \hbar^{-1} X^2 > 4 e^2 \pi / 3 =: L(0)$ (for $q=0$).

(b) For $q = 3$, such a steady state $\hhat$ exists if and only
if $\B \hbar^2 X^2 = E(0)$. For $q > 3$, $\hhat$ exists if and
only if $\B \hbar^{q-1} X^2 \leq E(0)$.  For all $q \geq 3$, if
$\hhat$ exists then $\E(\overline{h}) < \E(\hhat)$.
\end{theorem}

\begin{figure}[h] 
  \begin{minipage}{.5\linewidth} 
    \begin{center}
    \includegraphics[width=.5\textwidth]{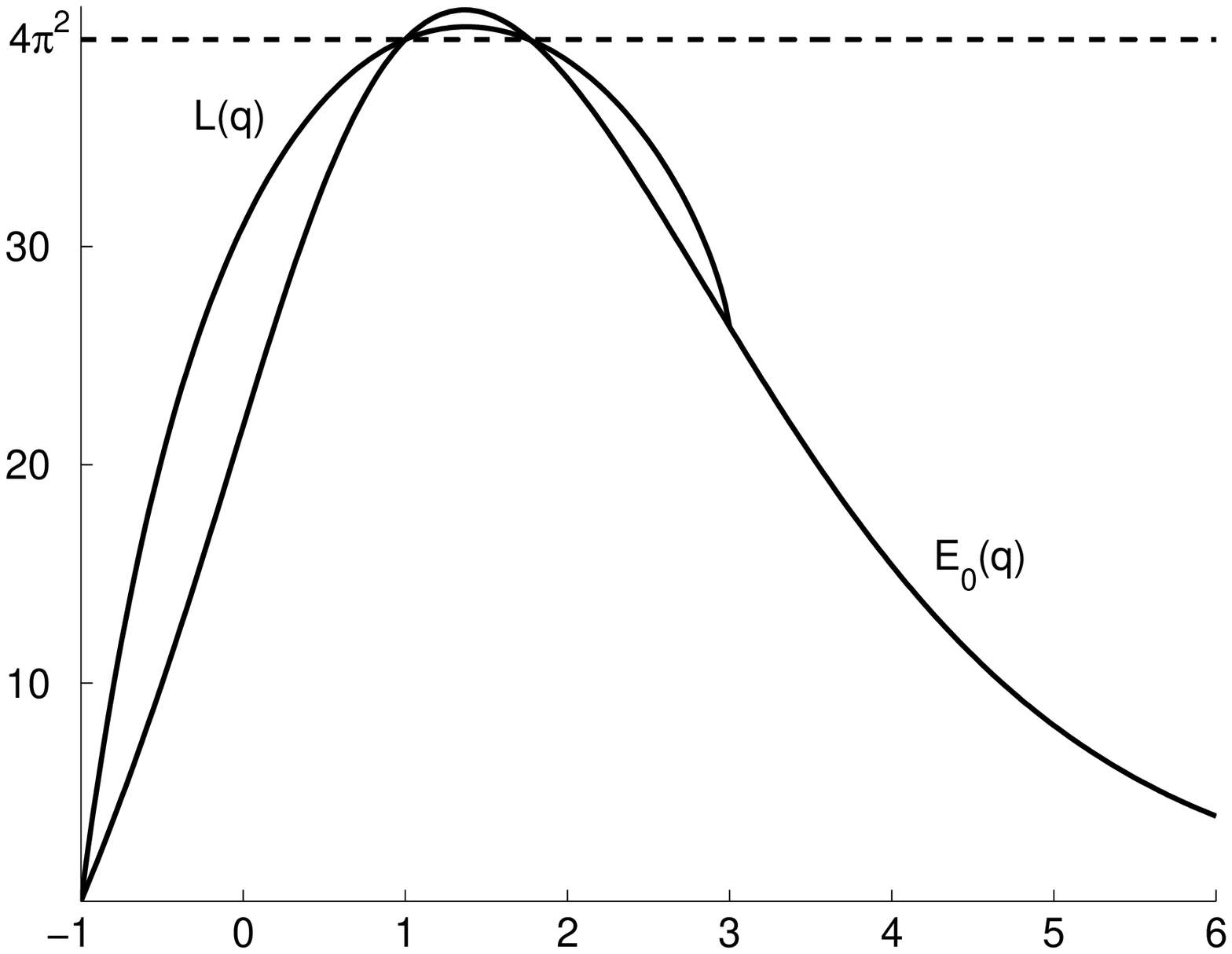}
  \vspace{-.2cm}  
  \caption{\label{constant_quantities} 
    Plots of $E_0(q)$ and $L(q)$.}
    \end{center}
  \end{minipage}%
  \begin{minipage}{.5\linewidth} 
    \begin{center}      
    \includegraphics[width=.5\textwidth]{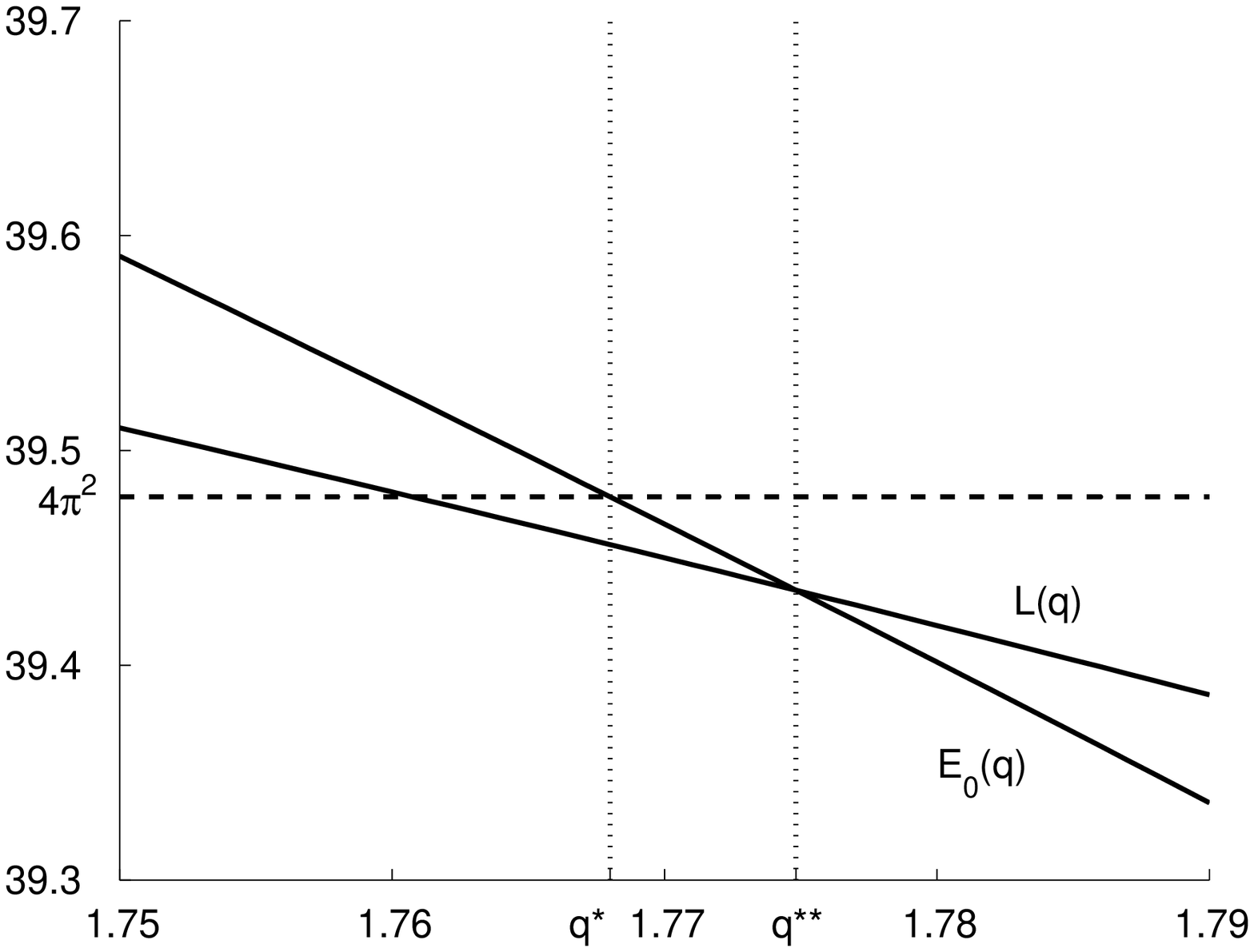}
  \vspace{-.2cm}  
  \caption{\label{closeup_quantities} Close-up view.} 
    \end{center}
  \end{minipage}
  \vspace{-0.2cm}
\end{figure}

To understand conditions (\ref{cd1}) and
(\ref{cd2}) see the plots of $E_0(q)$ and $L(q)$ in
Figures~\ref{constant_quantities} and \ref{closeup_quantities}
(constructed using the formulas for $A(0),P(0),E(0)$ in
\cite[\S3.1.2]{LP1}).  The graphs of $E_0$ and $L$ intersect at
$q=-1,1,3$ and at $q^{**} \approx 1.775$.  For $1<q<q^{**}$ the
figure suggests $E_0>L$, and so if the droplet steady state in Theorem~\ref{constant_droplet}(a) exists
then it certainly has lower energy than the constant steady
state.  On the other hand, it appears that $L > E_0$ when $-1<q<1$ and
when $q^{**} < q < 3$; in these cases the energy condition
(\ref{cd2}) may or may not hold when the existence condition
(\ref{cd1}) holds, so that the energy of the droplet steady
state (if it exists) might be higher or lower than that of the
constant steady state.

The dashed line at height $4\pi^2$ in
Figures~\ref{constant_quantities} and \ref{closeup_quantities}
intersects $E_0(q)$ at $1$ and $q^* \approx 1.768$, and intersects
$L(q)$ at $1$ and $1.761$ (approx.).  This line matters because the
constant steady state is a strict local minimum for the energy when
$\B \hbar^{q-1} X^2 < 4\pi^2$, by
Theorem~\ref{constant_stability}(b). For example, suppose that $2 \leq
q < 3$ and $\hbar$ and $X$ are such that $E_0(q)< L(q) < \B
\hbar^{q-1} X^2 < 4\pi^2$.  Then the constant steady state
$\overline{h}$ is a strict local minimum of the energy but is {\em
not} a global minimum since $\E(\hbar) > \E (\hhat )$ by
Theorem~\ref{constant_droplet}(a). A mountain pass idea then suggests
the energy might have a saddle point at which its value is greater
than $\E(\hbar)$. Such a saddle ought to be an energy unstable
positive periodic steady state, and should have period $X$ and area
$\hbar X$. In fact we already know such a positive periodic steady
state exists, by \cite[Theorem~12]{LP2} and
Theorem~\ref{periodic_constant} of this paper; we illustrate it in
Figure~\ref{q2.5_fig2}. Perturbing from the saddle in one direction
leads to relaxation to the constant (see the top part of
Figure~\ref{q2.5_fig2}), while perturbing in the opposite direction
yields apparent relaxation to a droplet (the bottom part of
Figure~\ref{q2.5_fig2}). The companion article \cite{LP4} contains
more extensive numerical investigations.

\section{Proofs of
Theorems~\protect\ref{nonlinear_unstable_general}--\ref{nonlinear_odd_power}
}
\label{proofs_nonlinear}
\subsection{Proof of Theorem~\protect\ref{nonlinear_unstable_general}}
\label{proof_nonlinear_unstable_general}
Recall that the positive periodic steady state $\h$ satisfies
equations (\ref{ss1}) and (\ref{ss2}):
\begin{equation} \label{ss6}
  \h^\ddd + r(\h) \h^\d = 0 , \qquad
  \h^\dd + H^\d(\h) = \text{const.},
\end{equation}
where $H^\dd(y) = r(y)$ is defined for $y>0$.

We compute the first four variations of the energy
(\ref{energy_is}). Let $u \in H^1(\T_X)$ have mean value zero. For
all $\varepsilon$ small enough that $\h + \varepsilon u > 0$
(so that $H(\h + \varepsilon u)$ makes sense), we have
\begin{eqnarray}
      \left. \frac{d\ }{d\varepsilon} \E(\h + \varepsilon u)              \right|_{\varepsilon = 0} 
& = & - \int_0^X \left[ \h^\dd + H^\d(\h) \right] u \, dx = 0 \label{var1} \\
&   &  \qquad \text{by (\ref{ss6}), since $u$ has mean value zero,} \nonumber \\
      \left. \frac{d^2\ }{d\varepsilon^2} \E(\h + \varepsilon u) 
                \right|_{\varepsilon = 0} 
& = & \int_0^X \left[ (u^\d)^2 - r(\h) u^2 \right] \, dx \label{var2} \\
& = & \text{numerator of Rayleigh quotient (\ref{rayleigh}) 
                                      for $\tau_1(\h)$,} \nonumber \\
      \left. \frac{d^3\ }{d\varepsilon^3} \E(\h + \varepsilon u) 
                \right|_{\varepsilon = 0} 
& = & -\int_0^X r^\d(\h) u^3 \, dx , \label{var3} \\
      \left. \frac{d^4\ }{d\varepsilon^4} \E(\h + \varepsilon u) 
                \right|_{\varepsilon = 0} 
& = & -\int_0^X r^\dd(\h) u^4 \, dx . \label{var4} 
\end{eqnarray}

If the steady state $\h$ is linearly unstable with respect to
zero--mean perturbations at period $X$, then the numerator of the
Rayleigh quotient (\ref{rayleigh}) is negative for some zero--mean
trial function $u \in H^1(\T_X) \setminus \{ 0 \}$.  We can assume $u
\in C^\infty(\T_X)$. From (\ref{var2}), the second variation of the
energy in the direction $u$ is negative, so that $\h$ is energy
unstable in the direction $u$ as desired.

From Theorems~1 and 3 in \cite{LP2}, 
$\h$ is linearly unstable with respect to zero--mean perturbations at
period $X$ if it is non-constant and: either the least period of $\h$
is $X/j$ for some integer $j \geq 2$, or $r=g/f$ is convex ($r^\dd
\geq 0$) and non-constant on the range of $\h$. Hence $\h$ is energy
unstable in those situations.

Now assume $\h$ is non-constant and $r$ is strongly convex ($r^\dd >
0$). We consider the second variation of $\E$ in the direction $u =
\pm \h^\dd$:
\begin{eqnarray}
\left. \frac{d^2\ }{d\varepsilon^2} \E(\h \pm \varepsilon
                \h^\dd) 
                \right|_{\varepsilon = 0} 
    & = & \int_0^X \left[ {\h^\ddd}^{\, 2}
          - r(\h) {\h^\dd}^{\! 2} \right] dx \nonumber \\
    & = & - \int_0^X r(\h) \left[ 
          \h^\d \h^\ddd 
          + {\h^\dd}^{\! 2} \right] dx 
          \qquad \mbox{by (\ref{ss6})} \nonumber \\
    & = & - \int_0^X r(\h) \left[ 
          \h^\d \h^\dd \right]^\d dx \nonumber 
    = \int_0^X r^\d(\h) 
          {\h^\d}^{\! \! 2} \h^\dd \, dx \\
    & = & \frac{1}{3} \int_0^X r^\d(\h) 
          \left[ {\h^\d}^{\! 3} \right]^\d dx 
    = - \frac{1}{3} \int_0^X 
          r^\dd(\h) {\h^\d}^{\! \! 4} \, dx . 
          \label{second_var}
\end{eqnarray}
By assumption, $r^\dd(\h(x)) > 0$ for all $x$.  Since $\h^\d$ is not
identically zero, the second variation in the direction $u = \pm
\h^\dd$ is negative, and so $\pm \h^\dd$ is an energy unstable
direction for $\h$.

It remains to prove $\h$ is energy unstable in the directions $u=\pm
\h^\d$. Here the second variation is zero, since (\ref{var2}) becomes
\[
\int_0^X 
\left[ (\h^\dd)^2 - r(\h) (\h^\d)^2 \right] \, dx 
 = -\int_0^X 
\left[ \h^\ddd + r(\h) \h^\d \right]
      \h^\d \, dx = 0 \qquad \text{by parts and (\ref{ss6})}.
\]

The third variation is zero for $u=\pm \h^\d$ by (\ref{var3}), because
$\h$ is even about its minimum point while $(\pm \h^\d)^3$ is odd, by
uniqueness for the nonlinear oscillator equation $\h^\dd + H^\d(\h) =
\text{const}$ in (\ref{ss6}). The fourth variation is negative by
(\ref{var4}), because $r^\dd(\h(x))>0$ by assumption.  Hence $\h$ is
energy unstable in the directions $u=\pm \h^\d$, completing the proof.
\qed

\vskip 6pt

A number of other authors, working on closely related topics, have
noted that translation invariance of the evolution implies the second
variation of the energy in the direction $\h^\d$ is zero
\cite{CGS,chafee,ANC90}. Those authors then proved energy instability
in the direction $u = \h^\d + \eta$ for some small function $\eta$.
Their arguments relied either on $\h$ having least period $X/j$ for
some $j \geq 2$ or else they did not impose a zero--mean requirement
on the perturbation. Above, we have used instead the assumption $r^\dd
> 0$.

\subsection{Proof of Theorem~\protect\ref{nonlinear_unstable_power}}
\label{proof_nonlinear_unstable_power}
\ 

When $q<1$ or $q>2$, Theorem~\ref{nonlinear_unstable_power} follows
from the last statement of Theorem~\ref{nonlinear_unstable_general}
since $r(y)=\B y^{q-1}$ is strongly convex.

When $q=2$ and $u=-\h^\dd$, formula (\ref{second_var})
shows the second variation of the energy is zero, since
$r(y)=\B y$ and $r^\dd \equiv 0$. The third
variation is negative because
\begin{eqnarray*}
\left. \frac{d^3\ }{d\varepsilon^3} \E(\h - \varepsilon \h^\dd) \right|_{\varepsilon = 0} 
& = & \B \int_0^X (\h^\dd)^3 \, dx \qquad \text{by (\ref{var3})} \\
& = & -2 \B \int_0^X \h^\d \h^\dd \h^\ddd \, dx \qquad \text{by parts} \\
& = & 2 \B^2 \int_0^X \h (\h^\d)^2 \h^\dd \, dx \qquad \text{since $\h^\ddd = - \B \h \h^\d$ by (\ref{ss6})} \\
& = & -\frac{2}{3} \B^2 \int_0^X (\h^\d)^4 \, dx \ < 0.
\end{eqnarray*}
Thus the steady state is energy unstable in the direction
$u=-\h^\dd$.  

Now suppose $q > 1$ and $E^\d(\alpha) > 0$. To obtain the unstable direction $u$, start by rescaling $k_\alpha$ to give
\begin{equation} \label{K_is}
K_\alpha(x) := \frac{P(\alpha)}{A(\alpha)} k_\alpha \left( P( \alpha)x \right). 
\end{equation}
By construction, $K_\alpha$ has period $1$ and mean value $1$. Then define
\begin{equation*} 
\kappa_\alpha(x) := \part{\;}{\alpha} K_\alpha(x);
\end{equation*}
$\kappa_\alpha$ is well-defined and smooth because $P$ and $A$ depend smoothly on $\alpha$ while $k_\alpha(x)$ is
jointly smooth in $(\alpha,x)$. Notice $\kappa_\alpha$ is even in $x$, has period $1$, and has mean value
zero:
\[
\int_0^1 \kappa_\alpha(x) \, dx 
= \part{\;}{\alpha} \int_0^1 K_\alpha(x) \, dx = \part{\;}{\alpha} (1) = 0.
\]
(See \cite[\S5.4]{LP2} for more properties of $\kappa_\alpha$.) Let $u(x) = \pm \kappa_\alpha(x/X)$. 

We want to show $u$ is an unstable direction. First we show that composing the rescalings of $\h$ to $k_\alpha$ and then $k_\alpha$ to $K_\alpha$ yields
\begin{equation} \label{rescaling_h_to_K}
\B X^2 \h(xX)^{q-1} = E K_\alpha(x)^{q-1} .
\end{equation}
Indeed, from the definition (\ref{K_is}) of $K_\alpha$, the righthand
side of (\ref{rescaling_h_to_K}) reduces to $P^2 k_\alpha(Px)^{q-1}$,
and then one can substitute for $k_\alpha$ in terms of $\h$ using
(\ref{rescaling}). Next one obtains the lefthand side of
(\ref{rescaling_h_to_K}) by using the relation $P = (\B/D)^{1/2q}
D^{1/2} X$ which relates the periods of $k_\alpha$ and $\h$
($P$ and $X$, respectively).

Now that we have (\ref{rescaling_h_to_K}), use (\ref{var2}) to compute the second variation of the energy in the direction $u(x)=\pm \kappa_\alpha(x/X)$ as 
\begin{eqnarray*}
\int_0^X \left[ \kappa_\alpha^\d(x/X)^2/X^2 - \B \h(x)^{q-1} \kappa_\alpha(x/X)^2 \right] dx
& = & 
\int_0^1 \left[ \kappa_\alpha^\d(x)^2 - \B X^2 \h(xX)^{q-1} \kappa_\alpha(x)^2 \right] dx/X \\
& = & 
\int_0^1 \left[ (\kappa_\alpha^\d)^2 - E K_\alpha^{q-1} \kappa_\alpha^2 \right] dx/X 
\qquad \text{by (\ref{rescaling_h_to_K})} \\
& < & 0 
\end{eqnarray*}
by the proof of \cite[Prop.~14]{LP2}, which uses both $q > 1$ and
$E^\d(\alpha) > 0$. Thus the steady state is energy unstable in the
direction $u=\pm \kappa_\alpha(x/X)$.  \qed

\subsection{Proof of Theorem~\protect\ref{nonlinear_stable_power}}
\label{proof_nonlinear_stable_power}
Since $E^\d(\alpha) < 0$ by hypothesis, Proposition~15 of \cite{LP2} implies $\mu_1(\alpha) \geq 0$, where
\begin{equation*} 
\mu_1(\alpha) := \min \left\{ 
\frac{ \int_0^1 \left[ (v^\d)^2 - E(\alpha) K_\alpha^{q-1} v^2 \right] dx}
       {\int_0^1 v^2 \, dx} 
    : v \in H^1(\T_1) \setminus \{ 0 \}, \int_0^1 v(x) \, dx = 0 \right\} .
\end{equation*}
(In fact $E^\d(\alpha) \leq 0$ would suffice to get this.)  Notice
$\mu_1(\alpha)=X^2 \tau_1(\h)$, by letting $v(x)=u(xX)$ and using the
identity (\ref{rescaling_h_to_K}) and the definition (\ref{rayleigh})
of $\tau_1$. Hence $\tau_1(\h) \geq 0$.

Consider $u \in H^1(\T_X) \setminus \{ 0 \}$ with mean value zero.
The first variation of $\E$ in the direction $u$ is zero by
(\ref{var1}), and the second variation of $\E$ in (\ref{var2}) is nonnegative because it equals the numerator of the Rayleigh quotient for $\tau_1(\h)$.

If the second variation is positive then $\h$ is energy stable in the
direction $\h$ and we are done. If the second variation is zero then
the Rayleigh quotient of $u$ is zero, and so $\tau_1(\h)=0$ and $u$
minimizes the Rayleigh quotient for $\tau_1$ in
(\ref{rayleigh}). Hence $u$ satisfies the Euler--Lagrange condition
$u^\dd + r(\h) u = \text{const}$, and so it satisfies $\LL u = 0$
where $\LL$ is the linearized operator defined in \cite[eq.~(4)]{LP2}
(take $a=0$ there).  Theorem~10(a) in \cite{LP2} and the hypothesis
$E^\d(\alpha) < 0$ now imply $u$ is a multiple of $\h^\d$, and so $u$
is odd.  Therefore the third variation (\ref{var3}) of $\E$ in the
direction $u$ is zero, because $u$ is odd and $\h$ is even. The fourth
variation is positive by (\ref{var4}) because $r(y)=\B y^{q-1}$ and
$r^\dd(y) < 0$ for $1<q<2$. Thus $\h$ is energy stable in the
direction $u$.  \qed

\subsection{Proof of Theorem~\protect\ref{nonlinear_odd_power}}
\label{proof_nonlinear_odd_power}
As usual, the first variation of $\E$ at $\h$ is zero by
(\ref{var1}).  We now prove non-negativity of the second variation, given by (\ref{var2}).

First note that $\h$ is symmetric about every point at which
$\h^\d = 0$, by uniqueness for the ODE $\h^\dd + H^\d(\h) = \text{const}$ (see (\ref{ss6})); the uniqueness theorem applies here since the
coefficient function $H^\d$ is $C^1$ (even $C^3$) on the range of the positive bounded function $\h$.  

Since $\h$ has a minimum at $x=0$ by hypothesis, we conclude $\h$ is
even and $\h^\d>0$ on $(0,X/2)$ (otherwise $\h$ would have period less than $X$).

Consider a minimizer of the Rayleigh quotient (\ref{rayleigh})
with respect to {\em odd} functions $u$; it is a smooth odd function
$\tilde{u}$ satisfying $\tilde{u}^\dd + r(\h)
\tilde{u} + \tilde{\tau} \tilde{u} = 0$ for some constant $\tilde{\tau}$.  Since
$\tilde{u}(0)=0$, one must have $\tilde{u}^\d(0) \neq 0$ because
otherwise $\tilde{u} \equiv 0$ by the uniqueness theorem for linear ODEs.
Also, $\tilde{u}(X/2)=0$ by the oddness and periodicity of $\tilde{u}$,
and so there is a point $b \in (0,X/2]$ with $\tilde{u}(b)=0$ and
$\tilde{u} \neq 0$ between $0$ and $b$.  Assume $\tilde{u} > 0$
between $0$ and $b$ (otherwise consider $-\tilde{u}$).  Then
\begin{eqnarray*}
  \tilde{\tau} \int_0^b \tilde{u} \h^\d \, dx 
  & = & - \int_0^b 
     \left[\tilde{u}^\dd + r(\h) \tilde{u} \right]
      \h^\d \, dx \qquad \text{since $\tilde{u}^\dd
        + r(\h) \tilde{u} + \tilde{\tau} \tilde{u} = 0$} \\ 
  & = & -\tilde{u}^\d(b) \h^\d(b) - \int_0^b
      \left[ \h^\ddd + r(\h) \h^\d \right]
     \tilde{u} \, dx \qquad \text{by parts, since $\h^\d(0)=0$} \\ 
  & = & -\tilde{u}^\d(b) \h^\d(b) \qquad \text{by (\ref{ss6})} \\
  & \geq & 0,
\end{eqnarray*}
because $\tilde{u}^\d(b) \leq 0$ and $\h^\d(b) \geq 0$.
Since $\tilde{u}$ and $\h^\d$ are positive on $(0,b)$ it
follows that $\tilde{\tau} \geq 0$. Hence the second variation of $\E$ in (\ref{var2}) is nonnegative, as desired.  

The third variation at $\varepsilon = 0$ is zero by (\ref{var3}), since $u$ is odd and $\h$ is even.

The fourth variation (\ref{var4}) is positive not just at $\varepsilon = 0$ but at all $\varepsilon \in (0,1)$:
\[
\frac{d^4\ }{d\varepsilon^4} \E(\h + \varepsilon u)
= -\int_0^X r^\dd(\h + \varepsilon u) u^4 \, dx > 0 
\]
by the strong concavity of $r$, provided $u \not \equiv 0$. Taylor's
theorem completes the proof, since for some $\tilde{\varepsilon} \in
(0,1)$,
\[
\E(\h + u) = \E(\h) + \frac{1}{2!}
\left. \frac{d^2\ }{d\varepsilon^2} \E(\h + \varepsilon u) 
\right|_{\varepsilon = 0} 
+  \frac{1}{4!} \left. \frac{d^4\ }{d\varepsilon^4} \E(\h + \varepsilon u) 
\right|_{\varepsilon = \tilde{\varepsilon}} > \E(\h).
\]
\qed

\section{Proofs of Theorems~\protect\ref{periodic_constant}---\protect\ref{constant_droplet}}
\label{proofs_energy}
\subsection{Proof of Theorem~\protect\ref{periodic_constant}}
\label{proof_periodic_constant}
We start by relating $\E(\h)$ to $\E(k_\alpha)$. Take 
\begin{equation} \label{H_is}
H(y) :=  \left\{ 
        \begin{array}{rl}
        \vspace{.05in}
        \ \frac{1}{q}\left[\frac{y^{q+1}}{q+1}-y\right], & q \neq 0, -1, \\
                \ y \log y - y,                         & q = 0, \\
                \ y - \log y,                   & q = -1 ,
        \end{array}
        \right.
\end{equation}
so that $H^\dd(y) = y^{q-1}$, and recall from the definition (\ref{energy_is}) of the energy that 
\begin{equation} \label{two_energies}
\E(\h) = \int_0^X \left[ \frac{1}{2} (\h^\d)^2 - \B H(\h) \right] dx \quad \text{and} \quad 
\E(k_\alpha) = \int_0^{P(\alpha)} \left[ \frac{1}{2} (k_\alpha^\d)^2 - H(k_\alpha) \right] dx .
\end{equation}

Denote the period of $\h$ by $\Pss=X$, and the area by $\Ass=\int_0^X \h \, dx$. Writing $P=P(\alpha)$ and $A=A(\alpha)$, the rescaling (\ref{rescaling}) implies
\renewcommand{\arraystretch}{1.5}
\begin{equation}  \label{scaling_PA}
P = 
        \left\{
                \begin{array}{cl}
                        \left( \frac{\B}{D} \right)^{\! \! 1/2q}
                        D^{1/2} \Pss ,        & q \neq 0, \\
                        e^{-D/2\B} \B^{1/2} \Pss , & q = 0,
                \end{array}
        \right.
\qquad \text{and} \qquad
A = 
        \left\{
                \begin{array}{cl}
                        \left( \frac{\B}{D} \right)^{\! \! 3/2q}
                        D^{1/2} \Ass ,        & q \neq 0, \\
                        e^{-3D/2\B} {\B}^{1/2} \Ass , & q = 0.
                \end{array}
        \right.
\end{equation}
\renewcommand{\arraystretch}{1}
Notice that the rescaling (\ref{rescaling}) can be written as
\begin{equation} \label{rescaling_rewritten}
\h(x) = \frac{\Ass}{A} \frac{P}{\Pss} \; k_\alpha \! \left( \frac{P}{\Pss} x \right) .
\end{equation}
From (\ref{scaling_PA}) we obtain the invariance relation
\begin{equation} \label{invariant}
\B \Pss^{3-q} \Ass^{q-1} = P^{3-q} A^{q-1}
= E(\alpha) ,
\end{equation}
and this implies
\begin{equation} \label{rescaling_bond}
\left(\frac{\Ass}{A}\right)^{\! \! 2}
\left(\frac{P}{\Pss}\right)^{\! \! 3}
= \B
\left(\frac{\Ass}{A}\right)^{\! \! q+1}
\left(\frac{P}{\Pss}\right)^{\! \! q} .
\end{equation}
Using (\ref{rescaling_rewritten}), (\ref{invariant}),
(\ref{rescaling_bond}) and the definitions (\ref{two_energies}), we at
last deduce a relation between $\E(\h)$ and $\E(k_\alpha)$:
\begin{equation} \label{energy_scaling}
\B^{-1} \overline{\h}^{\; -(q+1)} \frac{\E(\h) - \E(\overline{\h})}{\Pss}
= \overline{k_\alpha}^{\; -(q+1)} \frac{\E(k_\alpha) -
\E(\overline{k_\alpha})}{P(\alpha)} ,
\end{equation}
where the mean values are $\overline{\h} :=\Ass/\Pss$ and
$\overline{k_\alpha} :=A(\alpha)/P(\alpha)$. (When checking
(\ref{energy_scaling}), one can omit the linear terms in $H(y)$ from
the calculations, since $\h$ and $\overline{\h}$ have the same mean
value, as do $k_\alpha$ and $\overline{k_\alpha}$.)

In view of (\ref{energy_scaling}), then,
Theorem~\ref{periodic_constant} follows from:
\begin{proposition} \label{periodic_constant_prop}
Fix $\alpha_1 \in (0,1)$.

If $q \geq 2$ or $q < 1$ then $\E(k_\alpha) > \E(\overline{k_\alpha})$
for all $\alpha \in (0,1)$.

If $q = 1$ then $\E(k_\alpha) = \E(\overline{k_\alpha})$ for all
$\alpha \in (0,1)$.

If $1<q<2$ and $E^\d(\alpha) > 0$ $\forall \alpha \in (\alpha_1,1)$
then $\E(k_\alpha) > \E(\overline{k_\alpha})$ $\forall \alpha \in
[\alpha_1,1)$.

If $1<q<2$ and $E^\d(\alpha) < 0$ $\forall \alpha \in (\alpha_1,1)$
then $\E(k_\alpha) < \E(\overline{k_\alpha})$ $\forall \alpha \in
[\alpha_1,1)$.
\end{proposition}
Note that if $q \geq 2$ or $q < 1$ then $E^\d(\alpha)>0$ for all
$\alpha$ by \cite[Theorem~11]{LP2}.

\begin{proof}[Proof of Proposition~\protect\ref{periodic_constant_prop}]

The proof depends on a number of elementary differential equations and
inequalities that we derived in \S\S5.1--5.2 of \cite{LP2}, and the
reader may wish to skim those sections before proceeding.

If $q=1$ then $\E(k_\alpha) = \E(\overline{k_\alpha})$, as one sees
directly from the formula in (\ref{two_energies}), using that
$k_\alpha(x)=1+(\alpha-1) \cos x, P(\alpha)=2\pi$ and
$\overline{k_\alpha}=1$. So we assume $q \neq 1$ from now on, and
$\alpha \in (0,1)$.

By definition,
\begin{eqnarray}
  \frac{\E(k_\alpha) - \E(\overline{k_\alpha})}{P} & = & \frac{1}{P} \int_0^P \left[ \frac{1}{2} {k_\alpha^\d}^{\! 2} - H(k_\alpha) \right] dx  + H(A/P) \nonumber \\ 
& = & \frac{1}{P} \int_0^P {k_\alpha^\d}^{\! 2} \, dx - H(\alpha) + H(A/P) \qquad \text{by \cite[eq.~(21)]{LP2}.} \label{energy_diff_power}
\end{eqnarray}

First assume $q \neq -1$; then
\begin{equation} \label{energy_deriv_power}
\frac{d\ }{d\alpha} \frac{\E(k_\alpha) - \E(\overline{k_\alpha})}{P} = - \frac{P^\d(\alpha)}{P^2} \int_0^P {k_\alpha^\d}^{\! 2} \, dx + \left( 1 + (q+1) H(A/P) \left( \frac{A}{P} \right)^{\! \! -1} \right) \left( \frac{A}{P} \right)^{\! \! \d}\! \! (\alpha) ,
\end{equation}
where we have used \cite[eq.~(35)]{LP2} and the identity
$H^\d(y)=1 + (q+1)H(y) y^{-1}$ (valid for $q \neq -1$). Differentiating the function 
\begin{equation} \label{F_is}
\f(\alpha) := \left( \frac{A}{P} \right)^{\! \! -(q+1)} \frac{\E(k_\alpha) - \E(\overline{k_\alpha})}{P}
\end{equation}
(which is inspired by (\ref{energy_scaling})) with respect to
$\alpha$, we find from (\ref{energy_diff_power}) and
(\ref{energy_deriv_power}) that
\[
\f^\d(\alpha) = \left( \frac{A}{P} \right)^{\! \! -(q+2)} \! \left\{ \frac{1}{P} \int_0^P {k_\alpha^\d}^{\! 2} \, dx \cdot \left[ -(q+1) \left( \frac{A}{P} \right)^{\! \! \d} - \frac{A P^\d}{P^2} \right] + \left[ \frac{A}{P} + (q+1)H(\alpha) \right] \left( \frac{A}{P} \right)^{\! \! \d} \right\} .
\]
Substituting 
\[
\frac{A}{P} + (q+1)H(\alpha) = \frac{q+3}{2} \frac{1}{P} \int_0^P {k_\alpha^\d}^{\! 2} \, dx
\]
from \cite[eqs.~(31--32)]{LP2} yields
\begin{equation} \label{prior}
\f^\d(\alpha) = -\frac{1}{2} \left( \frac{A}{P} \right)^{\! \! -(q+2)} \frac{1}{P} \int_0^P {k_\alpha^\d}^{\! 2} \, dx \cdot \left[ (q-1) \left( \frac{A}{P} \right)^{\! \! \d} + 2 A P^{-2} P^\d \right] 
\qquad \mbox{for $q \neq -1$}.
\end{equation}

For $q=-1$ we obtain exactly the same formula (\ref{prior}) for
$\f^\d(\alpha)$, as follows:
\begin{eqnarray*}
      \f^\d(\alpha) 
& = & \frac{d\ }{d\alpha} \frac{\E(k_\alpha)- \E(\overline{k_\alpha})}{P} \\
& = & - \frac{P^\d}{P^2} \int_0^P {k_\alpha^\d}^{\! 2} \, dx + \left( 1 - P/A \right) \left( \frac{A}{P} \right)^{\! \! \d} \qquad \text{from (\ref{energy_diff_power}) and \cite[eq.~(35)]{LP2}} \\
& = & -\frac{1}{2} \left( \frac{A}{P} \right)^{\! \! -1} \frac{1}{P} \int_0^P {k_\alpha^\d}^{\! 2} \, dx \cdot \left[ -2 \left( \frac{A}{P} \right)^{\! \! \d} + 2 A P^{-2} P^\d \right]
\end{eqnarray*}
using that $(1-P/A)=\int_0^P (k_\alpha^\d)^2 \, dx/A$ when $q=-1$, by
\cite[eq.~(28)]{LP2}. The last equation is (\ref{prior}) for $q=-1$.

Equation (\ref{prior}) simplifies to
\begin{equation} \label{Fprime}
\f^\d(\alpha) = -\left[ \frac{1}{2} P^{-3} \left( \frac{A}{P} \right)^{\! \! -2q} \int_0^P {k_\alpha^\d}^{\! 2} \, dx \right] E^\d(\alpha) ,
\end{equation}
using that $E=P^2 (A/P)^{q-1}$ by (\ref{invariant}). Hence 
\begin{equation} \label{Eprime_Fprime} 
E^\d>0 \quad \Longleftrightarrow\quad \f^\d<0 \qquad \text{and} \qquad 
E^\d<0 \quad \Longleftrightarrow\quad \f^\d>0 .
\end{equation}
Also
\begin{equation} \label{limit}
\f(\alpha) \rightarrow 0 \quad \text{as} \quad \alpha \rightarrow 1 ,
\end{equation}
by the formula (\ref{energy_diff_power}) together with the
facts that $P,A \rightarrow 2\pi$ and $P/A \rightarrow 1$ as $\alpha
\rightarrow 1$ (see \cite[Lemma~6]{LP2}) and that $k_\alpha^\d \rightarrow 0$
uniformly as $\alpha \rightarrow 1$, by \cite[eq.~(21)]{LP2}.
Proposition~\ref{periodic_constant_prop} now follows from (\ref{Eprime_Fprime}), 
(\ref{limit}) and \cite[Theorem~11]{LP2} (which shows
$E^\d>0$ when $q \geq 2$ or $q < 1$).  For example, if $E^\d > 0$ on
$(\alpha_1,1)$ then $\f^\d < 0$ on $(\alpha_1,1)$; since
$\f(1) = 0$ we deduce $\f > 0$ on
$[\alpha_1,1)$, and so $\E(k_\alpha)-\E(\overline{k_\alpha})>0$ for $\alpha \in [\alpha_1,1)$.
\end{proof}
\subsection{Proof of Theorem~\protect\ref{periodic_droplet}}
\label{proof_periodic_droplet}
The proof relies on formulas for $A^\d(\alpha)$ and $E^\d(\alpha)$ that were derived in Lemmas~16 and 18 of \cite{LP2}: for all $q \neq -1$, 
\begin{equation} \label{A_prime_eq}
A^\d = 
-(q+1) H(\alpha) P^\d - \frac{q-1}{2} H^\d(\alpha) P ,
\end{equation}
\begin{equation} \label{E_prime_eq}
E^\d = - \left( \frac{A}{P} \right)^{\! q-2} 
  \left\{ P^\d \left[ (q-3) A + (q - 1) (q+1)H(\alpha) P \right]
        + \frac{1}{2} (q-1)^2 H^\d(\alpha) P^2 \right\}.
\end{equation}

Note that $E^\d>0$ if $q<1$ or $q \geq 2$, by \cite[Theorem~11]{LP2}.
Our assumptions therefore imply $-1 < q < 3$ and $E^\d(\alpha)>0$ for all $\alpha \in (0,1)$. Also $E(0)>0$ because $q>-1$ \cite[\S3.1.2]{LP1}. If we define 
\[
\hat{X} = \left[ E(0) / ( \B \Ass^{q-1} ) \right]^{1/(3-q)} 
\]
then
\[
\B \hat{X}^{3-q} \Ass^{q-1} = E(0) < E(\alpha) = \B X^{3-q} \Ass^{q-1}
\]
by (\ref{invariant}); here the value $\alpha \in (0,1)$ is determined by translating and rescaling $\h$ to a particular $k_\alpha$, as in \S\ref{power_law_review}. Hence $0 < \hat{X} < X$ (using that $q<3$). 

By rescaling the zero contact angle function $k_0$ on the interval $[0,P(0)]$ (as in \S\ref{power_law_review}; see \cite[Claim~5.1.2]{LP1} for details) we obtain a zero contact angle
droplet steady state $\hhat$ of $h_t = - (h^n h_{xxx})_x - \B(h^m h_x)_x$, with length $\hat{X}$ and area $\Ass$ as desired.

It remains to prove that the energy of this
droplet steady state is lower than the energy of the positive periodic
steady state $\h$. That is, we want to prove
\[
\int_0^X \left[ \frac{1}{2} (\hhat^\d)^2 - \B G(\hhat) \right] dx < \int_0^X \left[ \frac{1}{2} (\h^\d)^2 - \B G(\h) \right] dx 
\]
where
\begin{equation*} 
G(y) =  \left\{ 
        \begin{array}{rl}
        \vspace{.05in}
        \ \frac{y^{q+1}}{q(q+1)}, & q \neq 0, -1, \\
                \ y \log y - y,                         & q = 0, \\
                \ - \log y,                     & q = -1;
        \end{array}
        \right.
\end{equation*}
note that we can use $G$ instead of $H$ in the energy
because they differ only by a linear function (cf.~the definition (\ref{H_is}) of $H$) and $\hhat$ and $\h$ have the same area, $\Ass$. 

Since $\hhat$ is supported on $(0,\hat{X})$ and because $G(0)=0$,
the desired inequality is
\begin{equation} \label{desiderata}
\int_0^{\hat{X}} \left[ \frac{1}{2} (\hhat^\d)^2 - \B G(\hhat) \right] dx < \int_0^X \left[ \frac{1}{2} (\h^\d)^2 - \B G(\h) \right] dx .
\end{equation}
Next rescale $\hhat$ to $k_0$ and $\h$ to $k_\alpha$: from
(\ref{rescaling_rewritten}) and (\ref{rescaling_bond}) (with $\Pss$ replaced by $\hat{X}$ or $X$ as appropriate) we deduce (\ref{desiderata})
is equivalent to
\[
\left[ \frac{\Ass}{A(0)} \right]^{\! 2} \left[ \frac{P(0)}{\hat{X}} \right]^{\! 3} \int_0^{P(0)} \left[ \frac{1}{2} (k_0^\d)^2 - G(k_0) \right] dx
 <
\left[ \frac{\Ass}{A(\alpha)} \right]^{\! 2} \left[ \frac{P(\alpha)}{X} \right]^{\! 3} \int_0^{P(\alpha)} \left[ \frac{1}{2} (k_\alpha^\d)^2 - G(k_\alpha) \right] dx ,
\]
except that when $q=0$ we have to subtract
\[
\left[ \log \frac{\Ass}{A(0)} \frac{P(0)}{\hat{X}} \right] k_0 \quad \text{and} \quad \left[ \log \frac{\Ass}{A(\alpha)} \frac{P(\alpha)}{X} \right] k_\alpha
\]
from the integrands on the left and right sides, respectively. By substituting the relations 
\[
\hat{X} = \left[ E(0) / ( \B \Ass^{q-1} ) \right]^{1/(3-q)} \quad \text{and} \quad X = \left[ E(\alpha) / ( \B \Ass^{q-1} ) \right]^{1/(3-q)} 
\]
into the last inequality and using the definition $E(\alpha)=P(\alpha)^{3-q} A(\alpha)^{q-1}$, we see the desired inequality reduces to $\g(0) < \g(\alpha)$ where (for $q \neq 3$)
\begin{equation} \label{anotherG}
\g(\alpha) = 
A(\alpha)^{(q+3)/(q-3)} \int_0^{P(\alpha)} \left[ \frac{1}{2} (k_\alpha^\d)^2 - G(k_\alpha) \right] dx \quad \left( \text{$+ \frac{2}{3} \log A(\alpha)$, when $q=0$} \right).
\end{equation}
Thus to show the energy of $\hhat$ is lower than that of $\h$, it suffices to show $\g^\d(\alpha) > 0$ for all $\alpha \in (0,1)$, assuming  $-1 < q < 3$ and $E^\d(\alpha)>0$ for all $\alpha \in (0,1)$. 

To show $\g^\d > 0$, we substitute \cite[eq.~(32)]{LP2}
and \cite[eq.~(30)]{LP2} into the definition (\ref{anotherG}) of $\g$,
obtaining that
\[
\g = \frac{q-3}{q(q+3)} A^{2q/(q-3)} + \frac{q-1}{q+3} H(\alpha) P A^{(q+3)/(q-3)} \qquad \text{when $q \neq 3, 0, -1, -3$.}
\]
After differentiating this formula with respect to $\alpha$ and then
substituting for $A^\d(\alpha)$ from (\ref{A_prime_eq}), we
simplify with the help of (\ref{E_prime_eq}) to obtain
\begin{equation} \label{Gprime}
\g^\d = \frac{1}{q-3} H(\alpha) A^{6/(q-3)} \left( \frac{A}{P} \right)^{\! \! -q+2} E^\d \qquad \text{when $q \neq 3,0,-1,-3$.}
\end{equation}
When $q=0$ we obtain $\g = \frac{2}{3} - \frac{1}{3} H P A^{-1} +
\frac{2}{3} \log A$, by putting the $q=0$ versions of
\cite[eq.~(32)]{LP2} and \cite[eq.~(30)]{LP2} into
(\ref{anotherG}). Hence $\g^\d=- \frac{1}{3} H P^{-2} E^\d$, by
differentiating and using (\ref{A_prime_eq}) and
(\ref{E_prime_eq}). That is, (\ref{Gprime}) holds when $q=0$, also.

We conclude from (\ref{Gprime}) that $\g^\d > 0$ as desired (provided
$-1<q<3$ and $E^\d(\alpha)>0$ for all $\alpha \in (0,1)$), since for
$q>-1$ we have $H(\alpha)<0$ by the definition (\ref{H_is}). \qed
\subsection{Proof of Theorem~\protect\ref{nonzero_angle}}
\label{proof_nonzero_angle}
The proof involves rescaling arguments rather than the energy. Write
$\Pss=X$. Assume $\hhat$ is non-constant.

Suppose that in fact $\hhat$ is positive (and so smooth by
hypothesis). Then the least period of $\hhat$ equals $\Pss/j$
for some positive integer $j$, with the  area per period equaling
$\Ass/j$. If $j=1$ then $\hhat$ must be  a translate of
$\h$, by modifying slightly the uniqueness remarks in \cite[\S6.2]{LP2}
(using the assumption that $E^\d<0$ to get strict monotonicity of $E$). 

Thus we can assume $j \geq 2$. By rescaling $\h$ and $\hhat$ to $k_\alpha$ and $k_{\hat{\alpha}}$ for some $\alpha, \hat{\alpha} \in (0,1)$, as in (\ref{rescaling}), we get from (\ref{invariant}) that   
\[
E(\alpha) = \B \Pss^{3-q} \Ass^{q-1} \quad \text{and} \quad 
E(\hat{\alpha}) = \B (\Pss/j)^{3-q} (\Ass/j)^{q-1} .
\]
Hence
\begin{eqnarray*}
  4 \leq j^2 = \frac{E(\alpha)}{E(\hat{\alpha})} & < & \frac{E(0)}{E(1)}
  \qquad \text{since $E^\d < 0$} \\ & = & \frac{1}{4\pi^2}
  \frac{2}{q} (1+q) \; B \! \! \left(\frac{1}{2q}, \frac{1}{2}
\right)^{\! \! 3-q} B \! \! \left(\frac{3}{2q}, \frac{1}{2}\right)^{\!
  \! q-1} =: J(q) \ \text{say,} 
\end{eqnarray*}
by the formula for $E(0)$ in \cite[eq.~(3.13)]{LP1} and since $E(1)=P(1)^{3-q} A(1)^{q-1}=4\pi^2$ by \cite[Lemma~6]{LP2}. We will obtain a contradiction by showing $J(q)<4$, when $1<q \leq 1.75$; this will show $\hhat$ is not positive.

For $1 < q \leq 1.5$ we have
\[
J(q) \leq  
\frac{1}{4\pi^2} \frac{2}{1} (1+1.5) \; B \! \! \left(\frac{1}{2 \cdot 1.5}, \frac{1}{2} \right)^{\! \! 3-1}
B \! \! \left(\frac{3}{2 \cdot 1.5}, \frac{1}{2}\right)^{\! \! 1.5-1} \approx 3.17 < 4,
\]
where we have used that the Beta function $B(a,b) = \int_0^1 t^{a-1}
(1-t)^{b-1} \, dt$ is decreasing in its arguments, and is bigger than
$1$ when those arguments are less than $1$. For $1.5 < q \leq 1.75$ we
similarly have
\[
J(q) \leq  
\frac{1}{4\pi^2} \frac{2}{1.5} (1+1.75) \; B \! \! \left(\frac{1}{2 \cdot 1.75}, \frac{1}{2} \right)^{\! \! 3-1.5}
B \! \! \left(\frac{3}{2 \cdot 1.75}, \frac{1}{2}\right)^{\! \! 1.75-1} \approx 1.73 < 4,
\]
completing the contradiction. (Actually the best upper bound for $J(q)$ seems numerically to be about $1.04$, when $1 < q \leq 1.75$;
see \cite[Figure~6]{LP2}.)

The above contradiction implies $\hhat$ is not positive
everywhere. Consider therefore one component of the set $\{ x :
\hhat(x) > 0 \}$, say an interval with length $\hat{P} \leq
\Pss$. Write $\hat{A} \leq \Ass$ for the area under
$\hhat$ on this interval. Note the contact angles of
$\hhat$ must be the same at the two endpoints of the interval,
as a consequence of the nonlinear oscillator equation (\ref{ss3})
(see for example \cite[\S2.2]{LP1}). Suppose these contact angles
are zero, so that $\hhat$ rescales to $k_0$ on the interval,
using (\ref{rescaling}).  Then
\[
E(0) =   \B \hat{P}^{3-q} \hat{A}^{q-1} 
   \leq  \B \Pss^{3-q}   \Ass^{q-1} 
     =   E(\alpha)
\]
by applying (\ref{invariant}) twice, but this contradicts our
assumption that $E^\d < 0$. Thus the contact angles of $\hhat$ must
all be nonzero, as desired. \qed

\subsection{Proof of Theorem~\protect\ref{2_2_tango}}
\label{proof_2_2_tango}
Translate $h_{\rm ss1}$ and $h_{\rm ss2}$ so that they attain their
minimum values at $x=0$, and then rescale as in
\S\ref{power_law_review} to obtain $k_{\alpha_1}$ and $k_{\alpha_2}$
respectively.  Since $h_{\rm ss1}$ and $h_{\rm ss2}$ have the same
period and area, for which we write $\Pss=X$ and $\Ass$ respectively,
it follows from (\ref{invariant}) that
$E(\alpha_1)=E(\alpha_2)$. Notice $h_{\rm ss1} \neq h_{\rm ss2}
\Rightarrow \alpha_1 \neq \alpha_2$, in view of the expression
(\ref{rescaling_rewritten}) for $\h$ in terms of $k_\alpha$ and
$\Pss,\Ass,P(\alpha),A(\alpha)$. Since $E(\alpha_1)=E(\alpha_2)$ while
$E$ is strictly decreasing on $(0,\alpha_{crit})$ and strictly
increasing on $(\alpha_{crit},1)$, we conclude $\alpha_{crit}$ must
lie between $\alpha_1$ and $\alpha_2$.

We show $\alpha_1 < \alpha_2$.  The
hypothesis $\min h_{\rm ss1} < \min h_{\rm ss2}$ gives
\[
\alpha_1 \left( \frac{D_1}{\B} \right)^{\! \! 1/q} < \alpha_2 \left( \frac{D_2}{\B} \right)^{\! \! 1/q}
\]
by the rescaling (\ref{rescaling}). Next apply the first equation in
(\ref{scaling_PA}) to solve for $D_1$ in terms of $P(\alpha_1), \Pss$
and $\B$, and similarly solve for $D_2$ in terms of $P(\alpha_2),
\Pss$ and $\B$.  Substituting into the above inequality gives
$\alpha_1 P(\alpha_1)^{2/(q-1)} < \alpha_2 P(\alpha_2)^{2/(q-1)}$.
The strict increase of $\alpha \mapsto \alpha P(\alpha)^{2/(q-1)}$
implies $\alpha_1 < \alpha_2$.

Since $\alpha_1 < \alpha_{crit} < \alpha_2$, our hypothesis on $E^\d$
implies $E^\d(\alpha_1)<0$ and $E^\d(\alpha_2)>0$.
Theorem~\ref{nonlinear_stable_power} then implies that $h_{\rm ss1}$
is energy stable, and Theorem~\ref{nonlinear_unstable_power} implies
$h_{\rm ss2}$ is energy unstable.

Next we show $\E(h_{\rm ss1})<\E(h_{\rm ss2})$, or $\E(h_{\rm ss1}) - \E(\overline{h_{\rm ss1}}) <\E(h_{\rm ss2}) - \E(\overline{h_{\rm ss2}})$. In view of the
rescaling relation (\ref{energy_scaling}) for the energy, it suffices to prove
$\f(\alpha_1)<\f(\alpha_2)$, where $\f$ was
defined in (\ref{F_is}).

To prove $\f(\alpha_1)<\f(\alpha_2)$, we write
(\ref{Fprime}) as $\f^\d(\alpha) = \delta(\alpha)
(1/E)^\d(\alpha)$, where
\[
\delta(\alpha) = \frac{1}{2} P^3 A^{-2} \int_0^P {k_\alpha^\d}^{\! 2} \, dx .
\]
The point of this transformation is that $\delta(\alpha)$ is strictly
decreasing: $P^\d<0$ and $A^\d>0$, by \cite[Props.~7.3 and
7.4]{LP1}, while $\alpha \mapsto \int_0^P {k_\alpha^\d}^2 \, dx$
is decreasing by \cite[eq.~(35)]{LP2}. Also $(1/E)^\d>0$ on
$(0,\alpha_{crit})$ and $(1/E)^\d<0$ on $(\alpha_{crit},1)$, by
assumption. Hence
\[
\f^\d(\alpha) > \delta(\alpha_{crit}) (1/E)^\d(\alpha) \quad \text{for $\alpha \in (0,1), \alpha \neq \alpha_{crit}$.}
\]
Integrating this inequality from $\alpha_1$ to $\alpha_2$ and using
that $E(\alpha_1)=E(\alpha_2)$ gives $\f(\alpha_2)>\f(\alpha_1)$, as desired.

Finally, (\ref{Fprime}) shows $\f^\d(\alpha) < 0$ on
$(\alpha_{crit},1)$, and so $\f(\alpha_2)>\f(1)=0$ by
(\ref{limit}). Thus (\ref{energy_scaling}) yields $\E(h_{\rm ss2}) >
\E(\overline{h_{\rm ss2}})$. \qed

\subsection{Proof of Theorem~\protect\ref{constant_stability}}
\label{proof_constant_stability}
From the definition (\ref{rayleigh}) we see 
\[
\tau_1(\overline{h}) 
= \min_u \frac{\int_0^X (u^\d)^2 \, dx}{\int_0^X u^2 \, dx} - r(\hbar) = 
\left( \frac{2\pi}{X} \right)^{\! \! 2} - r(\hbar),  
\]
with the minimum being attained precisely for linear combinations of
$\sin(2\pi x/X)$ and $\cos(2\pi x/X)$. The first two paragraphs of the
theorem follow directly.

(a) Now suppose $r(\hbar) X^2 > 4\pi^2$. The variational formulas (\ref{var1}) and (\ref{var2}) in the proof of
Theorem~\ref{nonlinear_unstable_general} show the constant steady
state $\overline{h}$ is energy unstable in the sine and cosine directions,
since these are $\tau_1(\hbar)$-eigenfunctions and
$\tau_1(\hbar)<0$. Suppose next $r(\hbar) X^2 = 4\pi^2$ and
$r^\dd(\hbar)>0$. Then the first two variations of the
energy in the $\pm \sin$ directions are zero, by
(\ref{var1}) and (\ref{var2}). The third variation equals $r^\d(\overline{h})$ times the integral
of $\mp \sin^3$, by (\ref{var3}); thus the third
variation is also zero. The fourth variation is $- r^{\d
  \d}(\hbar) \int_0^X \sin^4(2\pi x/X) \, dx$, which is negative
because we assumed $r^\dd(\hbar)>0$. Thus the constant
steady state is energy unstable in the $\pm \sin$
directions. Argue similarly for the $\pm \cos$ directions.

(b) If $r(\hbar) X^2 < 4\pi^2$, or if $r(\hbar) X^2 = 4\pi^2$ and
$r^{\d \d}(\hbar)<0$, then $\tau_1(\hbar) \geq 0$ and so we get energy
stability by modifying the argument of part~(a) as follows. The first
variation of the energy in a direction $u$ is always zero. If the
second variation is positive then we are done. Otherwise the second
variation must be zero, so that $\tau_1(\hbar) = 0$ and $r^{\d
\d}(\hbar)<0$. Then $u$ must be a linear combination of sines and
cosines, and so the third variation is also zero. Then the fourth
variation is positive.

Furthermore, for $r(\hbar) X^2 < 4\pi^2$ we will prove $\hbar$ is a
strict local minimum of the energy, and is nonlinearly stable in
$H^1$. In doing this, we will use below a certain sufficiently small
number $\delta \in (0,1)$. Then considering $u \in H^1(\T_X)$ with
mean value zero and $\| u \|_{H^1(\T_X)}=1$, we find for all
$\varepsilon \in [0,\delta]$ that
\begin{eqnarray*}
\frac{d^2\ }{d\varepsilon^2} \E(\hbar + \varepsilon u)
& = & 
\int_0^X \left[ (u^\d)^2 - r(\hbar + \varepsilon u) u^2 \right] dx \\
& = & 
\delta \int_0^X \left[ (u^\d)^2 + u^2 \right] dx + (1-\delta) \int_0^X \left[ (u^\d)^2 - \frac{r(\hbar + \varepsilon u) + \delta}{1-\delta} u^2 \right] dx \\
& > & 
\delta \| u \|_{H^1(\T_X)}^2 + (1-\delta) \int_0^X \left[ (u^\d)^2 - \frac{4\pi^2}{X^2} u^2 \right] dx  
\geq \delta \| u \|_{H^1(\T_X)}^2 .
\end{eqnarray*}
In the second-to-last step above, we used $r(\hbar) X^2 < 4\pi^2$ and
$\|u\|_\infty \leq C \|u\|_{H^1}$ and we chose $\delta$ sufficiently
small (independent of $u$ and $\varepsilon$).

On the other hand, $|\varepsilon u| \leq \delta \| u \|_{L^\infty}
\leq \hbar/2$ provided $\delta$ is chosen small enough (independent of $u$), and hence
\begin{eqnarray*}
\frac{d^2\ }{d\varepsilon^2} \E(\hbar + \varepsilon u)
= 
\int_0^X \left[ (u^\d)^2 - r(\hbar + \varepsilon u) u^2 \right] dx 
& \leq & 
C(\hbar) \int_0^X \left[ (u^\d)^2 + u^2 \right] dx \\
&   & \text{for some constant $C(\hbar) \geq 1$} \\
& = & 
C(\hbar) \| u \|_{H^1(\T_X)}^2 .
\end{eqnarray*}

We deduce from the preceding estimates and the vanishing of the first
variation in (\ref{var1}) that if $u \in H^1(\T_X)$ has mean value
zero and $\| u \|_{H^1(\T_X)} \leq \delta$, then $\hbar + u > 0$ and
the energy varies quadratically away from $\hbar$:
\begin{equation} \label{quadratics}
\frac{1}{2} \delta \| u \|_{H^1(\T_X)}^2 \leq \E(\hbar + u) - \E(\hbar) \leq \frac{1}{2} C(\hbar) \| u \|_{H^1(\T_X)}^2 .
\end{equation}
The lefthand estimate implies $\hbar$ is a strict local
minimum of the energy, with respect to $X$-periodic zero-mean
perturbations.

We now prove $\hbar$ is 
nonlinearly stable, in the sense that if $h(x,t)$ is a smooth positive
solution of (\ref{evolve}) for $x \in \T_X$ and $t \in [0,T]$, for
some $T>0$, and if $h(\cdot,0)$ has mean value $\hbar$ and $\|
h(\cdot,0)-\hbar \|_{H^1(\T_X)} < \sqrt{\delta^3/4C(\hbar)}$, then $\|
h(\cdot,t)-\hbar \|_{H^1(\T_X)} < \delta/2$ for all $t \in
[0,T]$. Indeed, the quadratic bounds (\ref{quadratics}) and the
dissipation of the energy together imply
\begin{eqnarray*}
\frac{1}{2} \delta \| h(\cdot,t)-\hbar \|_{H^1(\T_X)}^2 & \leq & \E(h(\cdot,t)) - \E(\hbar) \\
& \leq & \E(h(\cdot,0)) - \E(\hbar) \leq \frac{1}{2} C(\hbar) \| h(\cdot,0)-\hbar \|_{H^1(\T_X)}^2 < \frac{1}{2} \frac{\delta^3}{4} ,
\end{eqnarray*}
so that $\| h(\cdot,t)-\hbar \|_{H^1(\T_X)} < \delta/2$ for all $t \in
[0,T]$. This stability result holds for all sufficiently small
$\delta$. \qed

\vskip 6pt This proof of nonlinear stability for a linearly stable
constant steady state does not carry over to linearly stable positive
periodic steady states $\h$, because there $\tau_1(\h)=0$ (due to
translational null directions). This zero eigenvalue is absent for the
constant steady state, since translation of $\hbar$ gives $\hbar$
again --- a trivial perturbation.  Imposing Neumann boundary
conditions eliminates the translational perturbations and their
associated zero eigenvalue.  Hence, the nonlinear stability proof {\em
would} hold for a positive steady state under the Neumann boundary
conditions $h_x=h_{xxx}=0$ at $x=0,X/2$ (see \S\ref{relation}),
provided the steady state is strictly linearly stable, {\it i.e.\ }the
first eigenvalue is positive.

\subsection{Proof of Theorem~\protect\ref{constant_droplet}}
\label{proof_constant_droplet}
Suppose first that $-1 < q < 3$. By (\ref{invariant}) with $\alpha=0$,
a steady state $\hhat$ supported on a single
interval of length $\hat{X} \leq X$ and with area
$\Ahatss=\hbar X$ and zero contact angles exists if and only if for
some length $\hat{X} \leq X$ we have 
\begin{equation} \label{invariant_0}
\B \hat{X}^{3-q} \Ahatss^{q-1} = E(0) = P(0)^{3-q} A(0)^{q-1} .
\end{equation}
That is, if and only if $\B X^{3-q} (\hbar X)^{q-1} \geq E(0)$,
which is (\ref{cd1}).

Suppose such a steady state $\hhat$ exists, supported say
on the interval $(0,\hat{X})$. By above,
\[
\hat{X} = 
\B^{1/(q-3)} P(0) [A(0)/\hbar X]^{(q-1)/(3-q)} .
\]
We want to show that $\E(\hhat)<\E(\overline{h})$ if and only if (\ref{cd2}) holds. For this, compute using $G$ like in \S\ref{proof_periodic_droplet} to find
\begin{eqnarray}
\E(\hhat) - \E(\overline{h})
& = & \int_0^X \left[ \frac{1}{2} \left( \hhat^\d \right)^2 - \B G(\hhat) + \B G(\hbar) \right] dx \nonumber \\
& = & \int_0^{\hat{X}} \left[ \frac{1}{2} \left( \hhat^\d \right)^2 - \B G(\hhat) \right] dx + \B G(\hbar) X \nonumber \\
& = & \frac{\Ahatss^2 P(0)^3}{A(0)^2 \hat{X}^3} \int_0^{P(0)} \left[ \frac{1}{2} (k_0^\d )^2 - G(k_0) \right] dx + \B G(\hbar) X \label{exception}
\end{eqnarray}
when $q > -1, q \neq 0$, by rescaling from $\hhat$ to $k_0$ and using (\ref{rescaling_rewritten})
and (\ref{rescaling_bond}) (with $\Pss$ replaced by $\hat{X}$, and $\Ass$ replaced by $\Ahatss$). Putting $\alpha=0$ into
\cite[eq.~(32)]{LP2} and \cite[eq.~(30)]{LP2} and using $H(0)=0$
enables us to evaluate $\int (k_0^\d)^2 \, dx$ and $\int G(k_0) \,
dx$, and hence we deduce
\begin{equation} \label{energy_difference}
\E(\hhat) - \E(\overline{h}) = \frac{\Ahatss^2 P(0)^3}{A(0)^2 \hat{X}^3} \frac{q-3}{q(q+3)} A(0) + \B G(\hbar) X \qquad \text{when $q>-1, q \neq 0$.}
\end{equation}
After substituting the definitions of $\hat{X}$ and
$\Ahatss=\overline{h} X$ from above, we find
\[
\E(\hhat) - \E(\overline{h}) = \frac{1}{q} \B^{3/(3-q)} (\hbar X)^{(3+q)/(3-q)} \left[ \frac{q-3}{q+3} A(0)^{2q/(q-3)} + \frac{1}{q+1} \left( \B \hbar^{q-1} X^2 \right)^{q/(q-3)} \right] 
\]
when $-1 < q < 3, q \neq 0$. Plainly now $\E(\hhat)<\E(\overline{h})$ if and
only if (\ref{cd2}) holds. 

When $q=0$ we find (\ref{exception}) has an
extra term $- \B \Ahatss \log \left[ \Ahatss P(0) / A(0)
\hat{X} \right]$, so that after again using
\cite[eq.~(32)]{LP2} and \cite[eq.~(30)]{LP2} and substituting for $\hat{X}$ and $\Ahatss$, we obtain
\[
\E(\hhat) - \E(\overline{h}) = \frac{1}{3} \B \hbar X \log \frac{A(0)^2 /e}{\B \hbar^{-1} X^2} \qquad \text{when $q = 0$.}
\]
Remembering that $A(0)=2 e^{3/2} \sqrt{\pi/3}$ when $q=0$, from
\cite[\S3.1.2]{LP1}, we conclude $\E(\hhat)<\E(\overline{h})$ if and
only if $\B \hbar^{-1} X^2 > 4 e^2 \pi / 3$. Incidentally, one can
check $L(q)$ is continuous at $q=0$.

For $q \geq 3$, simply modify the above proof from the case $-1<q<3$.
Notice when $q=3$ that (\ref{invariant_0}) becomes $\B (\hbar X)^2 =
E(0)$, which yields no formula for $\hat{X}$. And when $q \geq 3$ we
get $\E(\hhat) > \E(\overline{h})$ because the first term in
(\ref{energy_difference}) is nonnegative and the second is positive.
\qed

\section{Conclusions and Future Directions}
\label{conclusions}
If you perturb a positive periodic (or constant) steady state $\h$ of
the evolution equation (\ref{evolve}), without changing its area, then
towards which steady states might the solution subsequently evolve?

That is the broad question addressed by this paper. To answer it, we
focused on three specific questions:
\begin{itemize}
\item {\sc Existence:} Do there exist {\em other} steady states having
the same area and same period as $\h$, or having period a {\em
fraction} of the period of $\h$? If so, then these other steady states
are plausible contenders for the long time limit.  The constant steady
state $\overline{\h}$ obviously satisfies these requirements, but
there might be another positive periodic steady state (different from
$\h$) that does also, or perhaps a array of droplet steady states
having zero or nonzero contact angles.
\item {\sc Stability:} Are $\h$ or these other steady states linearly
stable? energy stable? If a steady state is the long time limit of
some generic solution, then surely it must be stable under
perturbations.
\item {\sc Relative energy levels:} Do any of these other steady
states have lower energy than $\h$? Only steady states with lower
energy are accessible, when starting from a small perturbation of
$\h$.
\end{itemize}
The {\sc existence} question was substantially answered for power law
coefficients by \cite[Theorem~12]{LP2}, Theorems~\ref{periodic_droplet} and
\ref{nonzero_angle} and \cite[Figures~3--5]{LP2},
also using Theorem~\ref{constant_droplet} when $\h$ is the
constant steady state. But the existence question remains open for
droplet steady states with {\em nonzero} contact angles, if we wish to
specify
the area and the length of the support. Some information on nonzero
angle droplet steady states is in our earlier paper \cite[\S5.2]{LP1}.

The {\sc stability} question was resolved for positive periodic steady
states in the power law case by \cite[Theorems~1,3,7,9]{LP2} and
Theorems~\ref{nonlinear_unstable_general}--\ref{nonlinear_odd_power}
here. In particular, Theorem~\ref{nonlinear_unstable_general} related
linear instability to energy
instability. Theorem~\ref{constant_stability} handled the case of
constant steady states. Our numerical simulations in the companion
article \cite{LP4} suggest that linearly unstable steady states are
indeed nonlinearly unstable, with the linear behavior dominating near
the steady state, but this observation is so far unsupported by a
`linearization theorem' for the power law evolution
$h_t = - (h^n h_{xxx})_x - \B(h^m h_x)_x$. (A linearization theorem {\em is} known in the
Cahn--Hilliard case $f \equiv 1$, by using semilinear operator
theory; see for example \cite[\S6]{ANC90}.)

The {\sc energy level} question has been largely settled in the
power law case by Theorems~\ref{periodic_constant}, \ref{periodic_droplet} and \ref{2_2_tango} when $\h$ is positive and periodic, and by
Theorems~\ref{periodic_constant} and \ref{constant_droplet} when $\h$ is
constant. When $\h$ has compact support with zero contact angle, use
Theorems~\ref{periodic_droplet} and \ref{constant_droplet}. For example,
Theorems~\ref{2_2_tango} suggests that when $m=n+0.77$, there can
exist two positive periodic steady states with the same period and
area, and that the unstable steady state has higher energy than the stable one. Our numerical simulations have
found heteroclinic connections from the high energy steady state to
the low energy one.

\vskip 6pt
\subsection*{Future directions.}
The stability question for {\em droplet} steady states (with zero and
nonzero contact angles) is open. So is the problem of computing
relative energy levels of non-zero angle droplet steady states with
respect to zero-angle droplets and constant and periodic steady
states.

Another open problem is to answer the existence, stability and energy
level questions for {\em general} coefficient functions $f$ and
$g$. We have treated power law coefficients, and Grinfeld and
Novick--Cohen \cite{GNC99} cover the Cahn--Hilliard equation, with $f
\equiv 1$ and $g(y)=1-3y^2$. But for general coefficients, about all
we know is that every positive periodic steady state is linearly and
energy unstable when $g/f$ is a convex function, by
Theorem~\ref{nonlinear_unstable_general}.

Finally, many of our existence, stability and relative energy level
theorems for the power law evolution would be improved if we knew
$E^\d(\alpha) < 0$ for all $\alpha$ when $1<q \leq 1.75$. We have not
been able to prove this conjecture, though numerically it is clear
from \cite[Figure~3]{LP2}.

\subsection*{{\bf Acknowledgments}}
The authors are grateful to Andrew Bernoff for stimulating comments on the energy landscape of phase space. 

Laugesen was partially supported by NSF grant number DMS-9970228, and
a grant from the University of Illinois Research Board. He is grateful
for the hospitality of the Department of Mathematics at Washington
University in St. Louis.

Pugh was partially supported by NSF grant number DMS-9971392, by the
MRSEC Program of the NSF under Award Number DMR-9808595, by the ASCI
Flash Center at the University of Chicago under DOE contract B341495,
and by an Alfred P. Sloan fellowship. Some of the computations were
done using a network of workstations paid for by an NSF SCREMS grant,
DMS-9872029.  Part of the research was conducted while enjoying the
hospitality of the Mathematics Department and the James Franck
Institute of the University of Chicago.


%
{\sc email contact:}  laugesen@math.uiuc.edu,  mpugh@math.upenn.edu
\end{document}